\documentclass[11pt]{article}
\usepackage{amsmath,amssymb}
\usepackage{pstricks}
\usepackage{pst-node}
\usepackage{graphicx}
\usepackage{pst-poly}
\usepackage{multido}
\usepackage{color} 

\newtheorem{result}{Theorem}
\newtheorem{deduce}{Corollary}

\newtheorem{define}{Definition}
\newtheorem{support}{Lemma}

\newtheorem{theorem}{Theorem}

\newcommand{\qed}{%
\ifmmode 
\else \leavevmode\unskip\penalty9999 \hbox{}\nobreak\hfill \fi
\quad\hbox{\qedsymbol}}
\newcommand{\openbox}{\leavevmode \hbox to.77778em{%
\hfil\vrule
\vbox to.675em{\hrule width.6em\vfil\hrule}%
\vrule\hfil}}
\newcommand{\qedsymbol}{\openbox}

\newcommand{\showgrid}{}
\newcommand{\gridon}{\renewcommand{\showgrid}{\psset{subgriddiv=1,griddots=10,gridlabels=6pt}\psgrid}}
\gridon

\begin{document}
\begin{center}
{\bf\LARGE  Revisiting the Hamiltonian Theme in the Square of a Block: The Case of $\mathbf{DT}$-Graphs }
\end{center}

\vskip8pt

\centerline{\large  Gek L. Chia$^{a,b}$  \ , \ Jan Ekstein$^{c}$ \ , \ Herbert Fleischner$^d$}

\begin{center}
\itshape\small$^a\/$Department of Mathematical and Actuarial Sciences, \\ Universiti Tunku Abdul Rahman,  Jalan Sungai Long, \\
Bandar Sungai Long, Cheras 43000 Kajang Selangor, Malaysia \\
$^b\/$ Institute of Mathematical Sciences, University of Malaya, \\ 50603 Kuala Lumpur,  Malaysia  \\ 
$^c\/$ Department of Mathematics, Institute for Theoretical Computer Science, and European Centre of Excellence NTIS - New Technologies
       for the Information Society \\ Faculty of Applied Sciences, University of West Bohemia, Pilsen, \\ Technick\'a 8, 306 14 Plze\v n, Czech Republic \\
$^d\/$ Institut f\"ur Computergraphik und Algorithmen 186/1, \\ Technical University of Vienna  \\ Favoritenstrasse 9–11, 1040 Wien, Austria \\
\end{center}

\vspace{3mm}

\begin{abstract}
The {\em square\/} of a graph $G\/$, denoted $G^2\/$, is the graph obtained from $G\/$ by joining by an edge any two nonadjacent vertices which have a common
neighbor. A graph $\/G\/$ is said to have the {\em ${\cal F} _k$ property\/} if for any set of $\/k\/$ distinct vertices $\{x_1, x_2, \ldots, x_k\}\/$  
in~$\/G\/$, there is a hamiltonian path from $x_1\/$ to $x_2\/$  in $\/G^2\/$ containing $k-2\/$ distinct edges of $G\/$ of the form 
$x_iz_i\/$, $i = 3, \ldots, k\/$. In \cite{fle3:refer}, it was proved that every $2\/$-connected graph has the ${\cal F}_3\/$ property. In the first part of 
this work, we extend this result by proving that every $2\/$-connected $DT\/$-graph has the ${\cal F}_4\/$ property (Theorem \ref{dt}) and will show in the 
second part that this generalization holds for arbitrary $2\/$-connected graphs, and  that there exist  $2\/$-connected graphs which do not have the 
${\cal F}_k\/$ property for any natural number $k \geq 5\/$. Altogether, this answers the second problem raised in \cite{cot:refer} in the affirmative.

\bigskip


\bigskip


\end{abstract} \vspace{3mm}

\section{Introduction and History}

All concepts not defined in this paper can be found in the book by Bondy and Murty, \cite{bm:refer}, or in the other references. However, we prefer 
definitions as given in Fleischner's papers if they differ from the ones given in \cite{bm:refer}. In particular, we  define a graph to be  {\em eulerian\/} 
if its vertices have even degree only; that is, it is not necessarily connected. This is in line with D.~K\"onig's original definition of an Eulerian graph, 
\cite{koe:refer}, and this is how eulerian graphs have been defined in Fleischner's papers quoted below (many authors call such graphs even graphs, whereas 
they consider a graph to be eulerian if it is a connected even graph). In any case, we consider finite loopless graphs only, but allow for multiple edges 
which may arise in certain constructions.

\vspace{1mm}
The study of hamiltonian cycles and hamiltonian paths in powers of graphs goes back to the late 1950s/early 1960s and was initiated by M. Sekanina who 
studied certain orderings of the vertices of a given graph. In fact, he showed in \cite{sek1:refer} that the vertices of a connected graph $G\/$ of order 
$n\/$ can be written as a sequence $a=a_1, a_2, \ldots, a_n=b\/$   for any given $a, b \in V(G)\/$, such that the distance $d_G(a_i,a_{i+1})  \leq 3\/$, 
$i=1, \ldots, n-1\/$. This led to the general definition of the {\em $k\/$-th power of a graph $G\/$\/}, denoted by $G^k\/$, as the graph with 
$V(G^k)=V(G)\/$ and $xy \in E(G^k)\/$ if and only if $d_G(x,y) \leq k\/$. Thus Sekanina's result says that $G^3\/$ is hamiltonian connected for every 
connected graph $G\/$.

\vspace{1mm}
Unfortunately, this result cannot be generalized to hold for $G^2\/$, the {\em square\/} of an arbitrary connected graph $G\/$ (the square of the subdivision 
graph of $K_{1,3}\/$ is not hamiltonian). Thus Sekanina asked in $1963\/$ at the Graph Theory Symposium in Smolenice, which graphs have a hamiltonian square, 
\cite{sek2:refer}. In $1964\/$, Neuman, \cite{neu:refer}, showed, however, that a tree has a hamiltonian square if and only if it is a caterpillar. On the 
other hand, it wasn't until $1978\/$ when it was shown in (\cite{und:refer}), that Sekanina's question was too general, for it was tantamount to asking 
which graphs are hamiltonian (that is, an $NP\/$-complete problem).

\vspace{1mm}
However, in $1966\/$ at the Graph Theory Colloquium in Tihany, Hungary, \mbox{ C. St. J. A.} Nash-Williams asked whether it is true that $G^2\/$ is 
hamiltonian if $G\/$ is $2\/$-connected, \cite{nas:refer}, and noted that L.W. Beineke and M.D. Plummer had thought of this problem independently as well.

\vspace{1mm}
By the end of $1970\/$, the third author of this paper answered Nash-Williams’ question in the affirmative; the corresponding papers \cite{{fle1:refer}, 
{fle2:refer}} were published in $1974\/$. In the same year, it was shown that this result implied that $G^2\/$ is hamiltonian connected for a $2\/$-connected 
graph $G\/$, \cite{chjkn:refer}.

\vspace{1mm}
Further related research was triggered by Bondy's question (asked in $1971\/$ at the Graph Theory Conference in Baton Rouge), whether hamiltonicity 
in $G^2\/$ implies that $G^2\/$  is {\em vertex pancyclic\/} (i.e., for every $v \in V(G)\/$ there are cycles of any length from $3\/$ through $|V(G)|\/$). 
In fact, Hobbs showed in $1976\/$, \cite{hob:refer}, that Bondy's question has an affirmative answer for the square of $2\/$-connected graphs and connected 
bridgeless $DT\/$-graphs (the latter type of graphs in which every edge is incident to a vertex of {\em d}egree {\em t}wo, was essential for answering Nash-
Williams' question – and it is essential for the main proofs of the current paper as well). The same issue of JCT B contains, however, a paper by Faudree and 
Schelp,  \cite{fs:refer},   in which they proved for the same classes of graphs, that since $G^2\/$ is hamiltonian connected, there are paths joining $v\/$ 
and $w\/$ of arbitrary length from $d_{G^2}(v,w)\/$ through $|V(G)| - 1\/$ for any $v, w \in V(G)\/$ (that is, $G^2\/$ is panconnected). They asked, however, 
whether this is a general phenomenon in the square of graphs (i.e., hamiltonian connectedness in $G^2\/$  implies panconnectedness in  $G^2\/$). Bondy's 
question and the question by Faudree and Schelp were answered in full in \cite{fle3:refer}.

\vspace{1mm}
Already in $1973\/$ (and published in $1975\/$) the most general block-cutvertex structure was determined such that every graph within this structure has 
a hamiltonian total graph, \cite{fh:refer}.

\vspace{1mm}
In the second part of the current work we establish in \cite{cf:refer} the strongest possible results in some sense (${\cal F}_k$-property), for the square 
of a block to be hamiltonian connected. As for hamiltonicity in the square of a block, the strongest possible result is cited Theorem \ref{fletheorem3} 
(\cite[Theorem 3]{fle3:refer}). Altogether, these results will enable us to establish (in joint work with others) the most general block-cutvertex structure 
such that if $G\/$ satisfies this structure then $G^2\/$ is hamiltonian connected or at least hamiltonian. That is, what has been achieved for total graphs, 
\cite{fh:refer}, will be achieved for general graphs correspondingly. Here, but also in the papers \cite{{fle1:refer},{fle2:refer}, {fle3:refer}, {fh:refer}}  
the concept of $EPS\/$-graphs plays a central role; and some of the theorems in the subsequent paper \cite{cf:refer}  require intricate proofs involving 
explicitly or implicitly $EPS\/$-graphs.

\vspace{1mm}
We are fully aware that there are shorter proofs on the existence of hamiltonian cycles in the square of a block; one has  been found by \v R{\'\i}ha, 
\cite{rih:refer}; and more recently, a still shorter proof was found by Georgakopoulos, \cite{geo:refer}. Moreover, a short proof of Theorem 
\ref{fletheorem3} (cited below) has been found by  M\"{u}ttel and Rautenbach, \cite{mr:refer}.       Unfortunately, their methods of proof do not seem to 
yield the special results which we can achieve with the help of $EPS\/$-graphs. This is not entirely surprising: \cite[Theorem 1]{fh:refer}  states that for 
a graph $G\/$, the total graph $T(G)\/$ is hamiltonian if and only if $G\/$ has an $EPS\/$-graph (note that the {\em total graph\/} of $G\/$ is the square of 
the subdivision graph of $G\/$).

\section{Preliminary Discussion}
By a {\em $uv\/$-path\/} we mean a path from $u\/$ to $v\/$. If a $uv\/$-path is hamiltonian, we call it a {\em $uv\/$-hamiltonian path\/}.

\begin{define}
Let $\/G\/$ be a graph and let $\/A = \{x_1, x_2, \ldots, x_k\}\/$ be a set of $\/k \geq 3\/$ distinct vertices  in $\/G\/$.  An   $x_1 x_2\/$-hamiltonian 
path in $\/G^2\/$ which contains $k-2\/$ distinct edges $x_iy_i \in E(G)\/$, $\/i = 3, \ldots, k\/$  is said to be $\/{\cal F} _k$.  Hence we speak of an  
${\cal F} _k$ $x_1x_2\/$-hamiltonian path. If $x_i\/$ is adjacent to $x_j\/$, we insist that $x_iy_i\/$ and $x_jy_j\/$ are distinct edges.  A graph $\/G\/$ 
is said to have the ${\cal F} _k$ property if for any set $\/A =\{x_1, \ldots, x_k\} \subseteq V(G)\/$, there is an  ${\cal F} _k$   $x_1x_2\/$-hamiltonian
path in $\/G^2\/$.
\end{define}

Let $G\/$ be a graph. By an {\em $EPS\/$-graph\/}, {\em $JEPS\/$-graph\/} respectively, of $G\/$, denoted $S = E \cup P\/$, $S = J \cup E \cup P\/$ 
respectively, we mean a spanning connected subgraph $S\/$ of $G\/$ which is the edge-disjoint union of an eulerian graph $E\/$ (which may be disconnected) 
and  a linear forest $P\/$, respectively a linear forest $P\/$ together with an open trail $J\/$.  For   $S = E \cup P\/$,  let $d_E(v)\/$ and  $d_P(v)\/$ 
denote the degree of $v\/$ in $E\/$ and  $P\/$, respectively. In the ensuing discussion we need, however, special types of $EPS$-graphs: thus 
a $[v;w]$-$EPS$-graph $S=E\cup P$ of $G$ with $v,w\in V(G)$, satisfies $d_P(v)=0$ and $d_P(w)\leq 1$. For $k\geq 2$, $[v;w_1,\dots,w_k]$-$EPS$-graphs are 
defined analogously, whereas in $[w_1\dots,w_k]$-$EPS$-graphs only $d_P(w_i)\leq 1$, $i=1,\dots,k$, needs to be satisfied. 

\vspace{2mm}
Let  $bc(G)\/$ denote the block-cutvertex graph of the graph $G\/$. If  $bc(G)\/$ is a path, we call $G\/$ a {\em block chain\/}. A block chain $G\/$ is 
called {\em trivial\/} if $E(bc(G)) = \emptyset\/$; otherwise it is called {\em non-trivial\/}. A block of $G\/$ is an endblock of $G\/$ if it contains at 
most one cutvertex of $G\/$.

\vspace{2mm}
In \cite[Lemma 2]{fle1:refer},  it was shown that  if $G\/$ is a block chain whose endblocks $B_1, B_2\/$ are $2\/$-connected and $v \in B_1\/$ and 
$w \in B_2\/$ are  not cutvertices  of $G\/$, then $G\/$ has an $EPS\/$-graph $S=E \cup P\/$ such that $d_P(v) =0=d_P(w)\/$. A more refined 
statement  is now given below. In Lemma 1 we apply \cite[Lemma 2, Theorem 3]{fle1:refer} and in Theorem 1 we apply Theorem~D (stated explicitly below) 
several times to the blocks of $G$, respectively to $G$ itself, to obtain $EPS$-graphs of the required type.

\vspace{2mm}
\begin{support}   \label{flelemma}
Suppose $G\/$ is a block chain with a cutvertex,  $v\/$ and $w \/$ are vertices in different endblocks of $G\/$ and are not cutvertices. Then

\vspace{1mm}
(i) there exists an $EPS\/$-graph $ E \cup P \subseteq G\/$ such that $d_P(v), \ d_P(w) \leq 1\/$. If the endblock which contains $v\/$ is 
$2\/$-connected, then we have $d_P(v) =0\/$ and $d_P(w) \leq 1\/$; and

\vspace{1mm}
(ii)  there exists a $JEPS\/$-graph $ J \cup E \cup P \subseteq G\/$ such that $d_P(v) =0 = d_P(w)\/$.  Moreover,  $v, w\/$ are the only odd vertices of 
$J\/$. Also, we have  $d_P(c) = 2\/$ for at most one cutvertex $c\/$ of $G\/$ (and hence $d_P(c') \leq 1\/$ for all other cutvertices $c'\/$ of $G\/$).
\end{support}

\vspace{2mm}
\noindent
{\bf Proof:} If $G\/$ is a path, the result is trivially true.

\vspace{1mm} 
So assume that $G\/$ is not a path. If $G\/$ has a suspended path (i.e., a maximal path whose internal vertices are 2-valent in $G$) starting at the 
endvertex $v\/$ of $G$, then let $P_v\/$ denote this path and let $v_1\/$ denote the other endvertex of $P_v\/$. Note that $v_1\/$ is a cutvertex 
of $G\/$. If there is no such suspended path, then define $ P_v\/$ to be an  empty path.  Likewise, $P_w\/$ is defined similarly with $w\/$ 
(respectively $w_1\/$) taking the place of $v\/$ (respectively $v_1\/$).

\vspace{1mm}
(i) By    \cite[Lemma 2]{fle1:refer}, $G'=G - ( P_v \cup  P_w)\/$ has an $EPS\/$-graph $S' = E' \cup P'\/$ with $d_{P'}(v_1)=0\/$ and $d_{P'}(w_1) \leq 1\/$. 
But this means that $G\/$ has an $EPS\/$-graph $S = E \cup P\/$ with $d_P(v) \leq 1\/$ and $d_P(w) \leq 1\/$ if we set $E = E'\/$ and 
$P = P' \cup P_v \cup P_w\/$. Clearly, in the case that $ P_v\/$ is an empty path, then $v=v_1\/$ and we have $d_P(v)=0\/$ and $d_P(w) \leq 1\/$.

\vspace{1mm}
(ii) Let $B\/$ be a block of $G\/$. Let $c_1,c_2\in V(B)$.  If $B\/$ is not an endblock, then let $c_1 , c_2 \in B\/$ be the  cutvertices of $G\/$ in $B\/$. 
If $B\/$ is an endblock of $G\/$, then let only one of $c_1, c_2\/$, say $c_2\/$, to be a cutvertex of $G\/$, and let $c_1=v$, $c_1=w$ respectively, 
depending on the endblock $c_1$ belongs to.   By  \cite[Theorem 3]{fle1:refer},  $B\/$ has a $JEPS\/$-graph $S_B = J_B \cup E_B \cup P_B\/$ with $d_{P_B}
(c_1) =0\/$,  $d_{P_B}(c_2) \leq 1\/$, and $c_1, c_2\/$ are the only odd vertices of $J_B\/$. If $B$ is not an endblock, then we may interchange $c_1$ and 
$c_2$. Thus we can ensure that for at most two blocks of $G$, $B'$ and $B''$ say, satisfying $B'\cap B''=c_2$, we have $d_{P_{B'}}(c_2)=d_{P_{B''}}(c_2)=1$. 

\vspace{1mm}
Note that if $B\/$ is not a $2\/$-connected block, then $E_B = \emptyset = P_B\/$ so that $S_B = J_B\/$. In this case, $d_{P_B}(c_1) = 0 = d_{P_B}(c_2)\/$.

\vspace{1mm}
By taking $S = \bigcup _B  S_B\/$, where the union is taken over all blocks $B\/$ of $G\/$, we have a $JEPS\/$-graph that satisfies the conclusion of (ii).

\vspace{1mm} This completes the proof.   \qed

\vspace{1mm}

\begin{result} \label{2edge}
Suppose $G\/$ is a $2\/$-connected graph and $v, w\/$ are two distinct vertices in $G\/$.
Then either

(i) there exists an $EPS\/$-graph $S=E \cup P \subseteq G\/$ with  $d_P(v)=0= d_P(w)\/$;

or

(ii) there exists a $JEPS\/$-graph $S= J\cup E \cup P \subseteq G\/$ with $v, w\/$ being the only odd vertices of $J\/$,  and   $d_P(v)=0= d_P(w)\/$.
\end{result}

\vspace{2mm}

\noindent
{\bf Proof:} If $G\/$ is a cycle, then clearly the result is true. Hence assume that $G\/$ is not a cycle.

\vspace{2mm}

Let $K'\/$ be a cycle in $G\/$ containing $v, w\/$. If $d_G(v)=2\/$, then we take a $[w;v]$-$EPS$-graph with $K' \subseteq E\/$. If $d_G(w) = 2\/$, then we 
take a $[v;w]$-$EPS$-graph with $K' \subseteq E\/$. In either case, Theorem D (stated below) guarantees the existence of such $EPS$-graphs. Thus conclusion 
(i) of the theorem is satisfied.

\vspace{2mm}
Hence we assume that $d_G(v), d_G(w) \geq 3\/$. We proceed by contradiction, letting $G\/$ be a counterexample with minimum  $|E(G)|\/$.

\vspace{2mm}
Let $G' = G-K'\/$ denote the graph obtained from $G\/$ by deleting all edges of $K'\/$ (including all possibly resulting isolated vertices).

\vspace{2mm}
(a) Suppose $G'\/$ is $2\/$-connected.  $G'\/$ either has an $EPS\/$-graph $S' = E' \cup P'\/$ or a $JEPS\/$-graph $S' = J' \cup E' \cup P'\/$ satisfying the 
additional property (i) or (ii), respectively.

\vspace{1mm} 
Suppose $S' = E' \cup P'\/$. Then set $E=K'\cup E'$, $P=P'\/$ to obtain an $EPS\/$-graph $S=E\cup P\/$ of $G\/$ satisfying property (i). If $G'$ has 
a $JEPS$-graph $S'=J'\cup E'\cup P'$ satisfying property (ii), then set $E=E'$, $P=P'$ and $J=J'\cup K'$, to obtain a $JEPS$-graph $S=J\cup E\cup P$ as
required. Whence $G'$ is not $2$-connected.

\vspace{2mm}
(b) Suppose $G'$ has an endblock $B'$ with $(B'-\gamma c')\cap \{v,w\}=\emptyset$ where $\gamma c'=c'$ if $B'$ contains a cutvertex $c'$ of $G'$, and 
$\gamma c'=\emptyset$ otherwise (in this latter case, $B'$ is a component of $G'$ having at least two vertices with $K'$ in common). It follows that 
$G'\supseteq H'$ where $H'$ is a block chain with $B'\subseteq H'$ and $G^*:=G-H'$ is $2$-connected. Suppose $H'$ is chosen in such a way that $G^*$ is as large 
as possible. 

\vspace{1mm}
It follows that if $H'$ is not $2$-connected then $|V(G^*)\cap V(H'-V(B'))|=1$. Denote the corresponding vertex with $c^*$ and observe that $c^*$ is 
a cutvertex if $c^*\in V(G')$. Also, by the choice of $B'$ and the maximality of $G^*$ we have 
$$(H'-c^*)\cap\{v,w\}=\emptyset$$
and $c^*$ is not a cutvertex of $H'$. Let $u'\in V(B')-\gamma c'$ be chosen arbitrarily. We set $\delta c^*=c^*$ if $c^*$ is a pendant vertex of $H'$, and 
$\delta c^*=\emptyset$ otherwise. By repeated application of Theorem \ref{thm1f} (see below) we obtain an $EPS$-graph $S'=E'\cup P'$ of $H'-\delta c^*\/$ 
with $d_{P'}(\delta c^*)=0$ (setting $d_{P'}(\emptyset)=0$) and $d_{P'}(u')\leq 1$. 

\vspace{1mm}
If however, $H'$ is 2-connected, i.e. $H'=B'$, then we let $c^{*}=(G'-B')\cap B'$, if $B'$ contains a cutvertex of $G'$, otherwise 
$c^{*}\in V(B')\cap V(K')$ arbitrarily. Futhermore we choose $u'\in V(B')-c^{*}$ arbitrarily. By Theorem \ref{thm1f}, $B'=H'$ has a $[c^{*};u']$-$EPS$-graph
$S'=E'\cup P'$.

\vspace{1mm}
Also, $G^*$ has an $EPS$-graph $S^*=E^*\cup P^*$ or a $JEPS$-graph 
$S^*=J^*\cup E^*\cup P^*$ with $d_{P^*}(v)=d_{P^*}(w)=0$; and $K'\subset E^*$, $K'\subset J^*\cup E^*$ respectively.

\vspace{2mm}
Observing that $P^*\cap P'=\emptyset$ and that $S^*$ and $S'$ are edge-disjoint, we conclude that $E=E^*\cup E'$ and $P=P^*\cup P'$ together with $J=J^*$ 
yield $S=E\cup P$, $S=J\cup E\cup P$ respectively, a spanning subgraph of $G$ as claimed by the theorem (observe that $d_P(c^*) = d_{P^*}(c^*)\/$ because 
$d_{P'}(c^*)=0\/$, and $d_{P^*}(c^*)=0\/$ if $c^* \in \{v, w\}\/$).

\vspace{2mm}
(c) Because of the cases solved already, we now show that $G'$ is connected and for every endblock $B'$ of $G'$, $V(B')\cap\{v,w\}\not =\emptyset$. For, if 
$G'$ is disconnected and because of case (b) already solved, $G'$ could be written as
$$G'=G'_1\ \dot \cup \ G'_2$$
 where $G'_i$ is a component of $G'$; and
$$G'_i\cap\{v,w\}\not =\emptyset,\quad i=1,2.$$
Without loss of generality $v\in G'_1$, $w\in G'_2$. Consequently, $G_i:=G'_i\cup K'$, $i=1,2$, is $2$-connected with $d_{G_1}(w)=2$, $d_{G_2}(v)=2$. Arguing 
as at the very beginning of the proof of this theorem (where we considered the case $d_G(v)=2$ or $d_G(w)=2$) we conclude that the corresponding $EPS$-graphs 
$S_i=E_i\cup P_i$ with $K'\subseteq E_i$, $i=1,2$, satisfy conclusion (i) of the theorem, and so does $S=E\cup P$ where $E=E_1\cup (E_2-K')$ and 
$P=P_1\cup P_2$.

\vspace{1mm}
Because of case (a) already solved, we thus have that $G'$ is a non-trivial block chain with $v,w$ belonging to different endblocks $B_v,B_w$ respectively, 
of $G'$ and they are not cutvertices of $G'$. Let $c_v$ and $c_w$ be the respective cutvertices of $B_v$ and $B_w$ (possibly $c_v=c_w$). If $B_v$ is not a 
bridge of $G'$ we use a $[v;c_v]$-$EPS$-graph $S_v$ of $B_v$ and a $[w;c_w]$-$EPS$-graph $S_w$ of $B_w$ if $B_w$ is also not a bridge, or $S_w=\emptyset$ if 
$B_w$ is a bridge. Proceeding similarly for every block $B$ of $G'-(B_v\cup B_w)$ we conclude that $G'$ has an $EPS$-graph $S'= E'\cup P'$ with 
$d_{P'}(v)=d_{P_v}(v)=0$ and $d_{P'}(w)=d_{P_w}(w)=0$, where $P_v\subseteq S_v$, $P_w\subseteq S_w$ (defining $d_{P_w}(w) =0\/$ if $P_w = \emptyset\/$). Thus 
in either case $S'\cup K'$ is an $EPS$-graph of $G$ satisfying conclusion (i). However, if both $B_v$ and $B_w$ are bridges, i.e., $d_{G'}(v)=d_{G'}(w)=1$, 
we introduce $z\notin V(G')$ and form $G_z:=G'\cup\{z,zv,zw\}$. $G_z$ contains a cycle $K_z$ through $z,v,w$ since $\kappa(G_z)\geq 2$, so it contains 
a $[v,w]$-$EPS$-graph $S_z=E_z\cup P_z$ with $K_z\subseteq E_z$. Trivially, $d_{P_z}(v)=d_{P_z}(w)=0$, and for the component $E_0\subseteq E_z$ with 
$z\in E_0$ we have $J:=(E_0-z)\cup K$ being an open trail joining $v$ and $w$. Setting $E=E_z-E_0$ and $P=P_z$ we conclude that $S=J\cup E\cup P$ is 
a $JEPS$-graph satisfying conclusion (ii) of the theorem. Theorem 1 now follows.
\qed

\vspace{3mm}
The following results from \cite{fh:refer}, \cite{fle1:refer}, and  \cite{fle3:refer}   will be used quite frequently in the proof of Theorem \ref{dt}.

\vspace{1mm} Let $G\/$ be a graph and let $W\/$ be a set of vertices in $G\/$. A cycle $K\/$ in $G\/$ is said to be  {\em $W$-maximal\/} if   
$|V(K') \cap W| \leq |V(K) \cap  W|\/$ for any cycle $K'\/$ of $G\/$. Moreover, we say that the $W$-maximal  {\em $K\/$ is $W$-sound\/} if 
$|V(K) \cap W| \geq 4\/$.

\vspace{1mm}
The following Theorems \ref{thm4fh} and \ref{thm3fh} are special cases of the theorems quoted.

\begin{theorem} (\cite[Theorem 4]{fh:refer}) \label{thm4fh}
Let $G\/$ be a $2\/$-connected graph and let $W\/$  be a set of five distinct vertices in $G\/$. Suppose $K\/$ is a $W\/$-sound cycle in $G\/$. Then there is 
an $EPS\/$-graph $S = E \cup P\/$ of $G\/$ such that $K \subseteq E\/$ and  $d_P(w) \leq 1\/$ for every $w \in W\/$.
\end{theorem}

An $EPS\/$-graph which satisfies the conclusion of Theorem \ref{thm4fh} is also called a {\em $W$-EPS-graph\/}.

\vspace{1mm}
\begin{theorem} (\cite[Theorem 3]{fh:refer}) \label{thm3fh}
Let $G\/$ be a $2\/$-connected graph and let $v, w_1, w_2, w_3\/$ be four distinct vertices of $G\/$. Suppose $K\/$ is a cycle in $G\/$ such that 
$\{v, w_1, w_2, w_3\} \subseteq K\/$. Then $G\/$ has a $[v; w_1, w_2, w_3]\/$-$EPS\/$-graph $S = E \cup P\/$ such that $K \subseteq E\/$.
\end{theorem}

\vspace{1mm}
Suppose $G\/$ is a $2\/$-connected graph and $v, w_1, w_2\/$ are distinct vertices in $G\/$. A cycle $K\/$ in $G\/$ is a {\em $[v; w_1, w_2]\/$-maximal cycle
\/} in $G\/$ if $\{v, w_1\} \subseteq V(K)\/$, and  $w_2 \in V(K)\/$ unless $G\/$ has no cycle containing all of $\{v, w_1, w_2\}\/$.

\begin{theorem} (\cite[Theorem 2]{fh:refer}) \label{thm2fh}
Let $G\/$ be a $2\/$-connected graph and let $v, w_1, w_2\/$ be three distinct vertices of $G\/$. Suppose $K\/$ is a $[v; w_1, w_2]\/$-maximal cycle 
in $G\/$. Then $G\/$ has a $[v; w_1, w_2]\/$-$EPS\/$-graph $S = E \cup P\/$ such that $K \subseteq E\/$.
\end{theorem}

\begin{theorem}   (\cite[Theorem 2]{fle1:refer}) \label{thm1f}
Let $G\/$ be a $2\/$-connected graph and let $v, w\/$ be two distinct vertices of $G\/$. Let $K\/$ be a cycle through $v, w\/$. Then $G\/$ has 
a $[v; w]\/$-$EPS\/$-graph $S= E \cup P\/$ with $K \subseteq E\/$.
\end{theorem}

\begin{theorem}  (\cite[Theorem 3]{fle3:refer}). Suppose $\/v\/$ and $\/w\/$ are two arbitrarily chosen vertices of a $\/2\/$-connected graph $\/G\/$.
Then $\/G^2\/$ contains a hamiltonian cycle $C\/$ such that the edges of $\/C\/$ incident to $\/v\/$ are in $\/G\/$ and at least one of
the edges of $\/C\/$ incident to $\/w\/$ is in $\/G\/$. Further, if $\/v\/$ and $\/w\/$ are adjacent in $\/G\/$, then these are three
different edges.  \label{fletheorem3}
\end{theorem}

\vspace{1mm}
A hamiltonian cycle in $G^2\/$ satisfying the conclusion of Theorem \ref{fletheorem3} is also called a $[v;w]\/$-hamiltonian cycle.  More generally, a 
hamiltonian cycle $C\/$ in $G^2\/$ which contains two edges of $G\/$ incident to $v\/$, and at least one  edge of $G\/$ incident to each 
$w_i\/$, $i=1, \ldots, k\/$,  is called a $[v; w_1, \ldots, w_k]\/$-hamiltonian cycle, provided the edges in question are all different.

\begin{theorem} (\cite[Theorem 4]{fle3:refer}). Let $G\/$ be a  $\/2\/$-connected graph. Then the following hold.

(i) $G\/$ has the ${\cal F} _3\/$ property.

(ii) For a given $q \in \{x, y\}\/$, $G^2\/$ has an $xy\/$-hamiltonian path containing an edge of $G\/$ incident to $q\/$.
     \label{fletheorem4}
\end{theorem}

\vspace{1mm}
By applying Theorems \ref{fletheorem3} and \ref{fletheorem4} to each block of a block chain $B\/$, we have the following.

\begin{deduce} \label{flecor}
Suppose $B\/$ is a  non-trivial  block chain with $|V(B)| \geq 3\/$ and  $v\/$ and $w \/$ are vertices in different  endblocks of $G\/$.  Assume further that
$v, w \/$ are not cutvertices of $B\/$. Then

(i) $B^2\/$ has a hamiltonian cycle which contains an edge of $B\/$ incident to $v\/$ and an edge of $B\/$ incident to  $w\/$. In the case that the endblock 
which contains $v\/$ is $2\/$-connected, then  $B^2$ has a hamiltonian cycle which contains two edges of $B$ incident to $v\/$ and an edge of $B\/$ incident 
to $w\/$.  Also,

(ii) $B^2\/$ has  a $vw\/$-hamiltonian path containing an edge of $B\/$ incident to $v\/$ and an edge of $B\/$ incident to $w\/$.
\end{deduce}

\section{$DT\/$-graphs}

Recall that a graph is called a {\em $DT\/$-graph\/} if every edge is incident to a $2\/$-valent vertex. If $G\/$ is a graph, we denote by $V_2(G)\/$ 
the set of all vertices of degree $2\/$ in $G\/$.

\vspace{2mm}
The following result which is interesting in itself, is obtained by applying Theorem~1 and the construction in \cite{fle1:refer} of a hamiltonian cycle/path 
in the corresponding spanning subgraph. 

\begin{deduce} \label{dt2edge}
Let $G\/$ be a $DT\/$-block and $x_1,x_2\in V(G)$ satisfying $N(x_1),  N(x_2) \subseteq V_2(G)\/$ and $x_1x_2\not \in E(G)$. Then  either (i) there exists 
a hamiltonian cycle in $G^2 -x_2\/$ whose edges incident to $x_1\/$ are in $G\/$, or else
(ii) there exists an $x_1x_2\/$-hamiltonian path in $G^2\/$ whose first and final edges are in $G\/$.
\end{deduce}

\vspace{1mm}
\begin{result}   \label{dt}
Every $\/2\/$-connected $DT\/$-graph  has the ${\cal F} _4$ property.
\end{result}

\vspace{3mm}
The proof of Theorem \ref{dt} is rather involved. We first give an outline of the general strategy used in the proof.

\vspace{1mm}
Let $G\/$ be a $2\/$-connected $DT\/$-graph and let $A = \{x_1, x_2, x_3 , x_4\}\/$ be a set of four distinct vertices in $G\/$.
Let $G^+\/$ denote the $2\/$-connected graph obtained from $G\/$ by adding a new vertex $y\/$ which joins $x_1\/$ and $x_2\/$. Then $G^+\/$ is a $DT\/$-graph 
unless $N_G(x_i) \not \subseteq V_2(G)\/$ for some $i \in \{1, 2\}\/$.  We shall show that $(G^+)^2\/$ contains a hamiltonian cycle $C\/$ containing edges of 
$G^+$ of the form $yx_1, yx_2, x_3z_3, x_4z_4\/$ where $x_3z_3, x_4z_4\/$   are edges of $G\/$. Then clearly $C\/$ gives rise to 
an ${\cal F} _4$ $x_1x_2\/$-hamiltonian path in $G^2\/$ when we delete the vertex $y\/$ from $(G^+)^2\/$.

\vspace{1mm}
In order to show the existence of such hamiltonian cycle $C\/$ in $(G^+)^2\/$, we shall apply induction or show that $G^+\/$ admits an $EPS\/$-graph 
$S = E \cup P\/$ with some additional properties.  In particular, in almost all cases, $E\/$ will contain a prescribed cycle $K^+\/$ passing through $y\/$. 
$K^+\/$ will also contain as many elements of $\{x_3, x_4\}\/$ as possible. Note that $G^+$ is $2$-connected and hence contains a cycle through $y$ and 
$x_i$, $i\in\{3,4\}$, which automatically contains $x_1,x_2$. 

\vspace{1mm}
Note that in \cite{fle1:refer} it was shown that if a $2\/$-connected $DT\/$-graph $H\/$ admits an $EPS\/$-graph, then $H^2\/$ has a hamiltonian cycle. We 
refer the reader to \cite{fle1:refer} for the method of constructing such hamiltonian cycle and to see how edges of $H\/$ can be included in such hamiltonian 
cycle.  Also, we may automatically assume that in an $EPS\/$-graph $S = E \cup P\/$ the edges of $P\/$ are the bridges of $S\/$ (otherwise, we could delete 
step-by-step $P\/$-edges (i.e., edges of $P\/$) until such situation is achieved).
 
\vspace{1mm}   
However, $G^+\/$ may not be a $DT\/$-graph and/or some elements in $A\/$ may be $2\/$-valent and (at least) one of its neighbors may not be $2\/$-valent. 
In such cases,  the existence of the various types of $EPS\/$-graphs $S\/$ in $G^+\/$ may not be sufficient to guarantee a hamiltonian cycle to begin with  
in $S^2\/$. Even if we can derive the existence of a hamiltonian cycle from these $EPS$-graphs, they may not suffice to guarantee a hamiltonian cycle  with 
the additional properties. Thus we need to consider neighbors of elements of $A\/$ to assure that they are incident to less than two $P\/$-edges. This 
applies, in particular, to $z_i \in N_{G}(x_i)\/$ with $z_{i}x_{i}\in E(K^{+})$, $i \in \{1,2,3,4\}\/$.

\vspace{1mm}
The following observations will be used quite frequently (sometimes implicitly) in the proof of Theorem \ref{dt}.

\vspace{5mm} \noindent
{\bf Observation (*):} {\em    Suppose $S= E \cup P\/$ is an $EPS\/$-graph of $G^+\/$ such that $d_P(x_i) \leq 1\/$ for $i=1, 2\/$.  Let $x\/$ be 
a $2\/$-valent vertex of $G\/$ belonging to $E\/$.

 \vspace{1mm}
 (i) Suppose $N(x) = \{u_{1}, u_{2}\}\/$. Then $S^2\/$ has a hamiltonian cycle which contains the edges $yx_1,yx_2$ and $u_{i}x\/$ for some $i\in \{1,2\}$   
     unless $x_j \in N(u_{j})\cup \{u_{j}\}\/$ and $d_P(x_j)=1, d_{S}(u_{j})>2$ for $j=1,2$; or for some $j\in \{1,2\}$, $d_P(x_j)=1, d_P(z_j)=2$ and 
     $z_{j}\in N(x_{j})\cap V(K^{+})$; or $d_P(u_1)=d_P(u_2)=2$ - in all three cases $N_G(x_j)\not\subseteq V_2(G)$.
 
 \vspace{1mm}
 (ii) We further note that any pendant edge in $S\/$ will always be contained in any hamiltonian cycle of $S^2\/$.

 \vspace{1mm}
 (iii) Consider $W\subseteq V(G^+)$ with $|W|=5$ and $K^+\subset G^+$. Suppose $|W\cap V(K^+)|\geq 4$. If $K^+$ is $W$-sound, then Theorem \ref{thm4fh} 
 applies. If, however, $K^+$ is not $W$-sound, then there is a $W$-sound cycle   $K^*\/$ with $W \subseteq K^*\/$ and we operate with $K^*\/$ in place of 
 $K^+\/$. This follows from the definition of $W\/$-soundness (see the discussion immediately preceding Theorem \ref{thm4fh}).                  }

\vspace{3mm}
The observations (i) and (ii) follow directly from the degree of freedom inherent in the construction of a hamiltonian cycle in $S^2\/$ as given 
in \cite{fle1:refer}.

\vspace{3mm}  The proof of Theorem \ref{dt} is divided into  several cases depending on whether $N(x_i)\subseteq V_2(G)\/$ or not, $i=1,2,3,4$.
Note that if $N(x_i) \not \subseteq V_2(G)\/$, then $d_G(x_i) =2\/$. If $d_G(x_i) = 2\/$,  we let $N(x_i) = \{u_i, v_i\}\/$ throughout the proof. Also, we 
define $x_i^* = x_i\/$ if $d_G(x_i) > 2\/$; and $x_i^*=z_{i}$ otherwise.

\vspace{2mm}
\begin{support}  \label{lemma4cd}
Let $G^+\/$ be defined as before with $N(x_3) \not \subseteq V_2(G)\/$ and $N(x_4) \not \subseteq V_2(G)\/$. Suppose $N(x_i) \subseteq V_2(G)\/$ for some 
$i \in \{ 1,2\}\/$. Assume further that every proper $2\/$-connected subgraph of $G\/$ has the ${\cal F} _4$ property. Then $(G^+)^2\/$ has a hamiltonian 
cycle containing the edges $x_1y, x_2y, x_3z_3, x_4z_4\/$ where $x_3z_3, x_4z_4\/$ are different edges of $G\/$.
\end{support}

\vspace{2mm}
\noindent
{\bf Proof:} \  By the hypotheses, $d_G(x_3)=d_G(x_4)=2$.  Assume without loss of generality that $N(x_1) \subseteq V_2(G)\/$.

\vspace{1mm}
{\bf (1)}  Suppose $\{u_i, v_i\} \neq \{x_1, x_2\}\/$ for $i = 3, 4\/$.

\vspace{1mm}
Let  $K^+\/$ be a cycle containing the vertices $y, x_1, x_2, x_4, u_4, v_4\/$.

\vspace{3mm} {\bf (1.1)} Assume that  $K^+\/$ also contains the vertex $x_3\/$.

\vspace{1mm}
We may assume that $$K^+ = yx_1z_1 \ldots u_4x_4v_4 \ldots u_3x_3v_3 \ldots z_2x_2y.$$

\vspace{1mm}  
{\bf (a)} Assume that  $u_4  \neq x_1\/$.

\vspace{2mm}
Since $\{x_1, x_2, x_3, x_4, u_3, u_4, z_2\} \subseteq V(K^+)\/$, Theorem \ref{thm3fh}  ensures the existence of a $[u_4; x_1,  z_2, x_2]\/$-$EPS\/$-graph 
$S_4 = E_4 \cup P_4\/$ of $G^{+}$ with $K^+ \subseteq E_4\/$ in the case $x_3x_4 \in E(G)\/$. Likewise, we obtain a $[u_4; x_1, u_3, x_2^*]$-$EPS\/$-graph 
$S_3 = E_3 \cup P_3\/$ of $G^{+}$ with $K^+ \subseteq E_3\/$ if $x_3x_4 \not \in E(G)\/$ where $x_2^* = x_2\/$ if $d_G(x_2) > 2\/$, and 
$x_2^* = z_{2} = V(K^+) \cap N_G(x_2)\/$ otherwise. It is straightforward to see that in both cases, the $EPS\/$-graph yields a hamiltonian cycle 
in $(G^+)^2\/$ as required by the lemma (see Observation (*)(i)).

\vspace{3mm}
{\bf (b)} Assume that  $u_4=x_1\/$ and $ v_3 = x_2\/$.

\vspace{2mm} \hspace{5mm} {\bf (b1)} Suppose $x_3\/$ and $x_4\/$ are adjacent or $N(x_3) \cap N(x_4) \neq \emptyset\/$.

\vspace{1mm}
(i) $x_3\/$ and $x_4\/$ are adjacent. Let $G^- = G - \{x_3, x_4\}\/$.   If $G^-\/$ is not $2\/$-connected, then it is a non-trivial block chain with 
$x_1, x_2\/$ belonging to different endblocks, and $x_1, x_2\/$ are not cutvertices of $G^-\/$. Hence  $(G^-)^2\/$ has a hamiltonian path $P(x_1, x_2)\/$ 
starting with an edge $x_1w_1\/$ of $G\/$ and  ending  with an edge $x_2w_2\/$ of $G\/$ (see Corollary \ref{flecor}(ii)). Then 
$$(P(x_1, x_2) - \{x_1w_1, x_2w_2\}) \cup \{x_1x_4, x_2x_3, x_4w_1, x_3w_2\}$$
defines a required ${\cal F}_4\/$ $x_1x_2\/$-hamiltonian path in $G^2\/$.

\vspace{1mm} 
If $G^-\/$ is $2\/$-connected, then $(G^-)^2\/$ has a hamiltonian cycle $C^-\/$ containing $x_1w_1, x_1t_1, \linebreak x_2w_2\/$ which are edges of $G\/$. 
Then $$(C^- - \{x_1w_1, x_1t_1,  x_2w_2\}) \cup \{w_1t_1, x_1x_4x_3w_2\}\/$$ is a required ${\cal F}_4\/$ $x_1x_2\/$-hamiltonian path in $G^2\/$.

(ii) Suppose $N(x_3) \cap N(x_4) = \{u\}\/$.

\vspace{1mm} 
If $d_G(u) = 2\/$, then let $G^- = G - \{x_3, x_4, u\}\/$ and proceed similarly as before to obtain a required ${\cal F}_4\/$ $x_1x_2\/$-hamiltonian path 
in $G^2\/$. Hence we assume that $d_G(u) > 2\/$. Suppose further that $G - x_i\/$ is   $2\/$-connected for some $i \in \{3, 4\}\/$. Then $G - x_i\/$ has 
the ${\cal F}_4\/$ property with $u\/$ taking the place of $x_i\/$; and any such ${\cal F}_4\/$ $x_1x_2\/$-hamiltonian path in $(G-x_i)^2\/$ can be extended 
to a required ${\cal F}_4\/$ $x_1x_2\/$-hamiltonian path in $G^2\/$. Thus we have to consider the case $\kappa(G-x_i) <2\/$ for $i \in \{3, 4\}\/$.

\vspace{1mm} 
Consider $G - x_4\/$. Since $d_G(x_4) =2\/$, $G' = G -x_4\/$ is a non-trivial block chain with $x_1, u\/$ belonging to different endblocks of $G'\/$ and are 
not cutvertices of $G'\/$. The endblock $B_u\/$ of $G'\/$ with $u \in V(B_u)\/$ also contains $x_3, x_2\/$ because $d_{G'}(u) \geq 2\/$ and 
$d_G(x_3) = d_{G'}(x_3) =2\/$. Hence $B_u\/$ is $2\/$-connected. Let $c\/$ be the cutvertex of $G'\/$ belonging to $B_u\/$.

\vspace{1mm} 
Suppose first $c \neq x_2\/$.   Because of the hypothesis of the lemma, $B_u\/$ has the ${\cal F}_4\/$ property. Correspondingly, there is a hamiltonian path 
$P(c, x_2)\/$ in $(B_u)^2\/$ containing $x_3w_3, uu'\/$ with $w_3 \in \{u, x_2\}\/$, which are different edges of $B_u\/$. Likewise, there is a hamiltonian 
path $P(x_1,c)\/$ in $(G'-B_u)^2\/$ by Theorem \ref{fletheorem4}, Corollary \ref{flecor}(ii), respectively. Then 
$$P(x_1, c) \cup (P(c, x_2) - uu') \cup \{u'x_4, x_4u\}\/$$ is a required ${\cal F}_4\/$ $x_1x_2\/$-hamiltonian path in $G^2\/$.

\vspace{1mm} 
Finally suppose $c= x_2\/$. By Theorem \ref{fletheorem4}(ii) or Corollary \ref{flecor}(ii), $(G' - B_u)^2\/$ has a $x_1x_2\/$-hamiltonian path $P_{1,2}\/$ 
ending with an edge $w_2x_2\/$ of $G\/$. By Theorem \ref{fletheorem3}, $(B_u)^2\/$ has a hamiltonian cycle $C_u\/$ with 
$\{ ux_3, x_2x_3, z_2x_2\} \subset E(B_u)\/$.
$$(P_{1,2} \cup C_u - \{w_2x_2, z_2x_2, ux_3\}) \cup \{w_2z_2, ux_4, x_4x_3\} $$ defines a hamiltonian path as required.

\vspace{2mm} 
\hspace{5mm} {\bf (b2)} Suppose $x_3\/$ and $x_4\/$ are not adjacent and $N(x_3) \cap N(x_4) = \emptyset\/$.

\vspace{1mm} 
Let $W = \{y, x_1, x_2, u_3, v_4 \}\/$. Then $K^+\/$ is $W\/$-sound. By Theorem \ref{thm4fh}, $G^+\/$ has an $EPS\/$-graph $S= E \cup P\/$ 
with $K^+ \subseteq E\/$ and  $d_P(w) \leq 1\/$ for every $w \in W\/$; and $d_P(x_3)=d_P(x_4)=0$. Because of the hypothesis of this case a required 
hamiltonian cycle can be constructed in $(G^+)^2\/$ (see Observation (*)(i)). In particular, the hamiltonian cycle contains $x_4v_4$ 
and $u_3x_3$.

\vspace{3mm}
{\bf (c)} Assume that   $u_4=x_1\/$ and  $v_3\ne x_2$.

\vspace{1mm}
If $x_3x_4 \not \in E(K^+)\/$, then Theorem  \ref{thm3fh} ensures the existence of an $[x_2; x_1, v_4, v_3]$-$EPS\/$-graph $S_3 = E_3 \cup P_3\/$ of $G^{+}$ 
with $K^+ \subseteq E_3\/$. By  construction, $S^2\/$ contains a hamiltonian cycle $C\/$ with $x_4v_4, x_3v_3 \in E(C)\/$. 
(see Observation (*)(i)). Hence we assume that $x_3x_4 \in E(K^+)\/$.

\vspace{1mm}
If $v_3x_2 \not \in E(K^+)\/$, or $v_3x_2 \in E(K^+)\/$ and $d_G(x_2) >2\/$, then we invoke Theorem \ref{thm2fh} to obtain 
a $[v_3; x_1, x_2^*]$-$EPS\/$-graph $S_3 = E_3 \cup P_3\/$ of $G^{+}$ with $K^+ \subseteq E_3\/$. If, however, $v_3x_2 \in E(K^+)\/$ and $d_G(x_2) = 2\/$, 
then Theorem \ref{thm2fh} ensures the existence of an $[x_2; x_1, v_3]$-$EPS\/$-graph $S_3 = E_3 \cup P_3\/$ of $G^{+}$ with $K^+ \subseteq E_3\/$. Note 
that $K^+\/$ contains all these special vertices. In all these cases, $(S_3)^2\/$ contains a hamiltonian cycle $C\/$ with $x_3x_4, x_3v_3 \in E(C)\/$ 
(see Observation (*)(i)).

\vspace{3mm} 
{\bf  (1.2)} In view of case (1.1), we may assume that $G^+\/$ has no cycle containing $y, x_4, x_3\/$, and that 
$$K^+ = yx_1z_1 \cdots u_4x_4v_4 \cdots z_2x_2y,$$ and $G^{+}-x_{3}$ is 2-connected if $G^{+}-x_{i}$ is 2-connected for some $i\in \{3,4\}$.

\vspace{2mm}
Without loss of generality, assume that $u_3 \not \in \{x_1, x_2\}\/$.

\vspace{3mm} 
{\bf (a)} Consider first the case that $G^* = G^+ - x_3\/$ is $2\/$-connected.

\vspace{1mm}
Define $W^* = \{ y, x_1, x_2^*, u_4, u_3\}\/$  if $x_1 \neq u_4\/$ and  $W^* = \{ y, x_1, x_2^*, v_4, u_3\}\/$  otherwise. Abbreviate 
$W^* = \{ y, x_1, x_2^*, t_4, u_3\}\/$ with $t_4 \in \{u_4, v_4\}\/$.

\vspace{1mm} 
\hspace{5mm} {\bf (a1)}  We first deal with the case $|W^*|=5\/$.

\vspace{1mm} 
In view of Observation (*)(iii), set $K^* = K^+\/$ if $K^+\/$ is $W^*\/$-sound in $G^*$, or else  there exists $K^* \supset W^* \/$ in $G^*$ 
(note $|K^+ \cap W^*| \geq 4\/$).

\vspace{2mm} 
\hspace{12mm} {\bf (a1.1)} Assume that $x_4 \in K^*$. In this case we may assume that $K^+=K^*$. By Theorem \ref{thm4fh}, there exists 
a $W^*$-$EPS$-graph $S^*=E^* \cup P^* $ of $G^*$ with $K^*\subseteq E^*$. Noting that  $d_{P^*}(u_3)\leq 1$, we set  $E=E^*$ and $P=P^* \cup \{u_3x_3\}$. 
Then  $S = E \cup P\/$ is an $EPS\/$-graph of $G^+\/$ whose structure implies that $(G^+)^2\/$ has a hamiltonian cycle containing the edges  $u_3x_3\/$ and 
$t_4x_{}\/$ (because $x_3\/$ is a pendant vertex in $S\/$ - see Observation (*)(i)-(ii)).

\vspace{2mm} 
\hspace{12mm} {\bf (a1.2)} Assume that $x_4 \not \in K^*\/$. Then $u_3\in K^*$ (hence $K^+\neq K^*$). Since $x_i\notin K^*$, for $i=3,4$, 
$d_G(t_4)>2$, $d_G(u_3)>2$. We define $x^{**}_2$ as $x^*_2$ with respect to $K^*$. 

First suppose $x^{**}_2=x^*_2$. By Theorem \ref{thm3fh}, there exists a $[u_4; x_1, u_3, x_2^*]$-$EPS\/$-graph $S^* = E^* \cup P^*\/$ of $G^*\/$ with 
$K^* \subseteq E^*\/$ if $x_1 \neq u_4\/$. By the same token, there is a $[u_4; v_4, u_3, x_2^*]$-$EPS\/$-graph $S^* = E^* \cup P^*\/$ of $G^*\/$ with 
$K^* \subseteq E^*\/$ if $x_1 = u_4\/$. In both cases, we set $E= E^*\/$, $P= P^*\cup \{x_3u_3\}\/$. Then $S= E \cup P\/$ is an $EPS\/$-graph of $G^+\/$ 
which yields a hamiltonian cycle in $(G^+)^2\/$ containing $u_3x_3\/$ and $x_4z\/$ for some $z \in N(x_4)\/$. If $x_4\/$ is a pendant vertex in $S^*\/$, 
then it is adjacent to $v_4\/$ (see Observation (*)(i)-(ii)).

If $x^{**}_2\neq x^*_2$, then we proceed analogously as before using $x^{**}_2$ instead of $x^*_2$. Note that $u_3=x^{**}_2$ is not an obstacle (we use 
Theorem \ref{thm2fh}) because of $d_S(x_2)=2$ since $d_G(x_2)=2$ and $x_2^*$ is also in $K^*$ (thus $x_2x_2^*\notin E(S)$) in this case.

\vspace{2mm} 
\hspace{5mm} {\bf (a2)} Assume  that $|W^*|= 4. \/$

\vspace{2mm}  
\hspace{12mm} {\bf (a2.1)} $W^{*}=\{y,x_{1},x_{2}^{*},t_{4}\}$ where $t_{4}\in \{u_{4},v_{4}\}$. If $u_{3}=t_{4}$, then we operate with 
a $[t_{4};x_{1},x_{2}^{*}]$-EPS graph $S^{*}=E^{*}\cup P^{*}$ of $G^*\/$ with $K^{+}\subseteq E^{*}$, which exists by Theorem C. If $u_{3}=x_{2}^{*}$, then 
we operate with a $[x_{2}^{*};x_{1},t_{4}]$-EPS graph $S^{*}=E^{*}\cup P^{*}$ of $G^*\/$ with $K^{+}\subseteq E^{*}$, which exists by Theorem C (note that 
$x_{2}^{*}\neq x_{2}$ in this case using $u_{3}\neq x_{2}$).

\vspace{1mm} 
In either case, set $E= E^*\/$ and $P = P^* \cup \{x_3u_3\}\/$. Then $S = E \cup P\/$ is an $EPS\/$-graph of $G^+\/$ which yields a hamiltonian cycle $C\/$ 
in $(G^+)^2\/$ containing either $x_3t_4x_4\/$ or $x_2^* x_3\/$,  $x_4t_4\/$ (see Observation (*)(i)).

\vspace{2mm} 
\hspace{12mm} {\bf (a2.2)}  $ W^* = \{ y, x_1, x_2^*, u_3 \}\/$. Then either (i) $u_4\neq  x_1\/$, or (ii) $u_4 = x_1\/$ and $v_4 \neq x_2^*\/$ or 
(iii) $u_4=   x_1\/$ and $v_4 = x_2^*\/$.

\vspace{1mm}  
In cases (i) and (ii) we are back to case (a2.1) with $u_{3}=t_{4}\neq x_2^*$.

\vspace{1mm} 
In case (iii) we have $x_2^* \neq x_2\/$ because $N(x_4) \neq \{x_1, x_2\}\/$. We consider $G' = G^+ - \{x_4, \delta x_2^*\}\/$; again, 
$\delta x_2^* = x_2^*\/$ if $x_2^*\/$ is a pendant vertex in $G^+ -x_4\/$ and $\delta x_2^* = \emptyset\/$ otherwise. Set $x_2'= x_2^*\/$ if 
$x_2^* \in V(G')\/$ and $x_2' = x_2\/$ otherwise. Suppose $\kappa(G') =1\/$. In any case, $G'\/$ has  different endblocks $B_1'\/$ and $B_2'\/$; they are 
$2\/$-connected with $x_1 \in B_1'\/$ and $x_2' \in B_2'\/$ not being cutvertices of $G'\/$. Since $G'\/$ is homeomorphic to $G\/$ if $x_2'= x_2\/$
(a contradiction to $\kappa(G') =1\/$), it follows that $x_2^* \in B_2'\/$ and that $x_2\/$ is a cutvertex of $G'\/$ since $\{x_2 \} = B_1' \cap B_2'\/$.
However, $3 = d_{G^+}(x_2) = d_{B_1'}(x_2) + d_{B_2'}(x_2) \geq 2+2\/$, an obvious contradiction. Thus $G'\/$ is $2\/$-connected in any case. Starting with
a cycle $K' \subseteq G'\/$ with $y, x_1, x_2, x_3 \in V(K')\/$ we apply Theorem \ref{thm2fh} to obtain an $[x_1; u_3, x_{2}^{**}]$-$EPS\/$-graph 
$S' = E' \cup P'\/$ of $G'\/$ with $K' \subseteq E'\/$, where $x_{2}^{**}=N_{G}(x_{2})-x_{2}^{*}$. Setting 
$E = E'\/$, $P = P' \cup \{x_1x_4, \delta(x_{4}x_{2}^{*})\}\/$, where $\delta(x_{4}x_{2}^{*})=x_{4}x_{2}^{*}$ if $x_{2}^{*}\notin V(G')$ and 
$\delta(x_{4}x_{2}^{*})=\emptyset$ otherwise, we obtain $S= E \cup P\/$ of $G^{+}$ with $K' \subseteq E\/$ and $d_P(x_1) =1\/$ and $d_P(u_3) \leq 1\/$. 
It is clear that $S^2\/$ yields a hamiltonian cycle of $(G^+)^2\/$ as required (see Observation (*)(i)).

\vspace{2mm} 
\hspace{5mm} {\bf (a3)} Assume  that $|W^*|= 3. \/$

\vspace{2mm} 
Then  $W^* =\{ y, x_1, x_2^*\}\/$.

\vspace{1mm} 
Hence $u_{3}\notin \{x_{1},y\}$, therefore $u_{3}=x_{2}^{*}$. Analogously $t_{4}\notin \{x_{1},y\}$, therefore $t_{4}=v_{4}=x_{2}^{*}$.
That is, $u_{3}=x_{2}^{*}=t_{4}=v_{4}$ and $x_{1}=u_{4}$. $G' = G^+ - x_4\/$ is $2\/$-connected since  there is a cycle $K'\/$ in $G'\/$ containing $y\/$ 
and $x_3\/$ and hence also $x_2^*, x_1, v_3\/$. If $v_3 \neq x_1\/$, we operate with an $[x_2^*; x_1, v_3]$-$EPS\/$-graph $S' = E' \cup P'\/$  of $G'\/$ 
with $K' \subseteq E^*\/$ (by Theorem~\ref{thm2fh}). Setting $E = E^*\/$ and $P = P^* \cup \{x_4x_2^*\}\/$, we obtain an $EPS\/$-graph 
$S = E \cup P\/$ of $G^{+}$ which will yield a hamiltonian cycle in $(G^+)^2\/$ containing  $x_3v_3, x_4v_4\/$ (see Observation (*)(i)). If $v_3 = x_1\/$, 
then $G-x_4\/$ is $2\/$-connected (since $N(x_3) = N(x_4)\/$). Hence  $G - x_4\/$ has the ${\cal F}_4\/$ property with $v_4\/$ taking the place of $x_4\/$; 
and any such ${\cal F}_4\/$ $x_1x_2\/$-hamiltonian path in $(G-x_4)^2\/$ can be extended to a required ${\cal F}_4\/$ $x_1x_2\/$-hamiltonian path in $G^2\/$.
This finishes the proof of case {\bf (a)}.  

\vspace{3mm}
{\bf  (b)} \
Now consider the case where $G^* = G^+ - x_3 \/$ has a cutvertex and hence $G^+-x_4$ has also a cuvertex, because of the assumptions of case \textbf{(1.2)}. 
Thus $G^*$ is a non-trivial block chain since $d_G(x_3)=2$. Note that $K^+\/$ is contained in some endblock $B_y\/$ of $G^*\/$. 

\vspace{1mm} 
Let $W = \{y, x_1, x_2^*, x_3, t_4\}\/$ where we define $t_4$ as follows: 
\begin{itemize}
 \item $t_4=u_4$ if $u_4\neq  x_1$;
 \item $t_4=v_4$ if $u_4=x_1$ and either $x_2^*=x_2$ or $v_4\neq x_2^*\neq x_2$;
 \item $t_4=x_2$ if $u_4=x_1$ and $v_4=x_2^*\neq x_{2}$.
\end{itemize}
Note that by this definition of $t_4$, $|W|=5$. 

\vspace{1mm} 
Assume first that the cycle $K^+\/$ (which passes through  $y, x_1, x_2, x_2^*, u_4, x_4, v_4\/$) is $W$-sound in $G^+$. Let $\widehat{G}\/$ denote 
the subgraph of $G^+\/$ which is a non-trivial block chain containing $u_3, x_3, v_3\/$ such that $G^+ - \widehat{G} = B_y \/$. Suppose $w_3\/$ is 
the vertex in one of the endblocks of $\widehat{G}$ and $w'_3$ the vertex in the other endblock of $\widehat{G}$ such that 
$\widehat{G}\cap B_y=\{w_3, w'_3\}$. Possibly $\{w_3, w'_3\} \cap \{u_3, v_3\} \neq  \emptyset\/$, but $\{w_3, w'_3\} \neq \{u_3, v_3\}\/$.


\vspace{1mm} 
We replace $\widehat{G}\/$ in $G^+\/$ by a path $P_4 = a_1a_2x_3a_3a_4\/$ (where $a_1, a_3\/$ are identified  with $w_3, w'_3\/$ respectively, and 
$\{a_2,a_3\}=\{u_3,v_3\}$) to obtain the graph  $G''\/$. Note that $K^+\subseteq G''\/$. Set $W = \{y, x_1, x_2^*, x_3, t_4\}\/$ as above. Then $K^+\/$ is 
$W\/$-sound  (by assumption), and by Theorem \ref{thm4fh}, $G''\/$ has an $EPS\/$-graph $S'' = E'' \cup P''\/$ such that $K^+ \subseteq E''\/$ and 
$d_{P''}(z) \leq 1\/$ for every $z \in W\/$.

\vspace{1mm} 
(\textbf{b1)} Suppose $E(P_4) \cap E(P'') = \emptyset\/$. Then $P_4 \subseteq E'' \/$. Since $\widehat{G}\/$ is a non-trivial block chain, 
by Lemma \ref{flelemma}(ii), $\widehat{G}\/$ contains a $JEPS\/$-graph $\widehat{S} = \widehat{J} \cup \widehat{E} \cup \widehat{P}\/$ such that 
$d_{\widehat{P}}(w_3)=0=d_{\widehat{P}}(w'_3)$, and $w_3,w'_3$ are the odd vertices of $\hat J$; hence $d_{\widehat{J}}(x_3)=2$ and $d_{\widehat{P}}(x_3)=0$. 
Note that by the second part of Lemma \ref{flelemma}(ii) we can make sure that $\mbox{min}\{d_{\widehat{P}}(u_3),d_{\widehat{P}}(v_3)\}\leq 1$. In this case, 
we obtain an $EPS\/$-graph $S = E \cup P\/$ of $G^+\/$ by setting $E = (E'' -P_4) \cup \widehat{J} \cup \widehat{E} \/$ and $P = P'' \cup \widehat{P}\/$. 
Here $d_{P}(x_3) =0\/$ and $d_{P}(w) \leq 1\/$ for every $w \in W -x_3\/$.

\vspace{1mm} 
\textbf{(b2)} Suppose $E(P_4)\cap E(P'')\neq\emptyset$. That is, $V(P_4)\subseteq V(P'')$ (so that $E(P_4) \cap E(E'')= \emptyset$) and $d_{P''}(x_3) =1\/$. 
This means that either $a_2x_3 \not \in E(P'')\/$ or $x_3a_3 \not \in E(P'')\/$. Suppose $x_3a_3 \not \in E(P'')\/$ (so that $a_3a_4 \in E(P'')\/$). 
In this case, we delete $x_3v_3$ from $\widehat{G}$ and thus split $\widehat{G}$ into two block chains $\widehat{G_1}$ and $\widehat{G_2}$ with 
$x_3,w_3 \in \widehat{G_1}$ and $v_3,w'_3\in \widehat{G_2}$. If $\widehat{G_j}$ is an edge only, then $\widehat{S_j}=\widehat{G_j}$. If 
$\widehat{G_2}=\emptyset$, then $\widehat{S_2}=\emptyset$. Otherwise by Lemma~\ref{flelemma}(i) (or by Theorem \ref{thm1f} if $\widehat{G_2}$ is 
2-connected), $\widehat{G_j}$ has an $EPS$-graph $\widehat{S_j}=\widehat{E_j} \cup \widehat{P_j}$ where $d_{\widehat{P_1}}(w_3) \leq 1\/$, 
$d_{\widehat{P_1}}(x_3) = 1$, $d_{\widehat{P_2}}(v_3) \leq 1$, $d_{\widehat{P_2}}(w'_3)\leq 1$, $j=1,2$. Now, if we take 
$E=\widehat{E_1}\cup\widehat{E_2} \cup E''\/$ and $P = \widehat{P_1} \cup \widehat{P_2} \cup (P'' -\{a_2,a_3\})\/$, we have an $EPS\/$-graph $S=E \cup P$ 
of $G^+\/$ with $d_{P} (w) \leq 1\/$ for every $w \in W\/$ (note that $w_3a_2, w'_3a_3 \in P'')$, $x_3$ is a pendant vertex in $S$, and it works also if 
$\widehat{G}$ is a path on at least 4 vertices.

In both cases \textbf{(b1)} and \textbf{(b2)}, a required hamiltonian cycle in $(G^+)^2\/$ can be constructed from $S$ (see Observation (*)(i)-(ii)).
Note that $G^+-x_4$ is 2-connected if $K^+=yx_1x_4x_2^*x_2y$ (hence $d_{G}(x_2)=2$) and if $d_{G}(x_2^*)>2$. Here we have a contradiction to the assumption 
of this case \textbf{(1.2)(b)}.

\vspace{2mm} 
Now assume that the cycle $K^+$ is not $W$-sound. Since $y,x_1,x_2^*,t_4\in K^+$ and $|W|=5$, there exists a cycle $K^*\subseteq G^+$ containing all of $W$
and not $x_4$.

\vspace{1mm} 
(i) Suppose $t_4=v_4$ or $t_4=x_2$. In both cases, $K^*$ contains $u_4=x_1$ and $v_4=t_4$, $v_4=x_2^*$, respectively, but not $x_4$. Hence $G^+-x_4$ is 
2-connected, a contradiction with assumptions of case \textbf{(1.2)(b)}.

\vspace{1mm} 
(ii) Suppose $t_4=u_4$. Because $G^+-x_3$ has a cutvertex, without loss of generality suppose that $u_3\notin K^+$ but clearly $u_3\in K^*$. Hence 
$u_3\notin \{x_1,x_2^*,u_4\}$. We define $x^{**}_2$ as $x^*_2$ with respect to $K^*$.

First suppose $x^{**}_2=x^*_2$. By Theorem \ref{thm3fh}, $G^+$ has a $[u_4; x_1, x_2^*,u_3]$-$EPS$-graph $S=E\cup P$. Note that either $x_4$ is a pendant 
vertex in $S$, or else $x_4$ is a vertex in $E$. It is clear that $S^2$ yields a hamiltonian cycle of $(G^+)^2$ as required (see Observation (*)(i)-(ii)).

If $x^{**}_2\neq x^*_2$, then we proceed analogously as before using $x^{**}_2$ instead of $x^*_2$. Note that $x^{**}_2=u_3$ or $x^{**}_2=u_4$ is not 
an obstacle (we use Theorem \ref{thm2fh}) because of $d_S(x_2)=2$ since $d_G(x_2)=2$ and $x_2^*$ is also in $K^*$ (thus $x_2x_2^*\notin E(S)$); and 
$x^{**}_2=u_3=u_4$ is not possible in this case.

\vspace{3mm}
{\bf (2)} \ Suppose  $\{u_3, v_3\} = \{x_1, x_2\}\/$.

\vspace{1mm} 
Note that, in $G\/$, there exists a cycle containing $x_3\/$ and $x_4\/$ (and hence also the vertices $u_3, v_3, u_4, v_4\/$).

\vspace{1mm}
Let $G^* = G^+ - x_3\/$ which is homeomorphic to $G\/$ and thus  $G^*\/$ is $2\/$-connected. Note that there exists a cycle $K^*=K^+\/$ (see above) in 
$G^*\/$ containing the vertices $y, x_1, x_2, x_4\/$.

\vspace{3mm} 
{\bf (2.1)} Suppose $w \in N(x_4) - \{x_1, x_2\}\/$ exists; let $x_2z_2 \in E(K^*)\/$. Note that $d_{G^*}(x_2) =2\/$ if $d_G(z_2) >2\/$.

\vspace{1mm}
By Theorem \ref{thm3fh}, there exists an $[x_1; x_2, z_2, w]\/$-$EPS\/$-graph, an $[x_1; x_2, z_2]\/$-$EPS\/$-graph by Theorem C respectively, if $w=z_2$;
in both cases we denote $S^* = E^* \cup P^* \subset G^*\/$ with $K^* \subseteq E^*\/$ and $d_{P^*}(x_1) =0\/$. Note that $K^*$ is $[x_1;x_2,z_2]$-maximal 
if $z_2=w$. Set $E = E^*\/$ and $P = P^* \cup \{x_1x_3\}\/$; thus $d_P(x_1)=1\/$. Also, $d_P(x_2)+d_P(z_2)\leq 1$ since $d_P(z_2)>0$ implies $d_P(x_2)=0$
since $d_{G^*}(x_2)=2$. Then $S = E \cup P\/$ is an $EPS\/$-graph of $G^+\/$ and a hamiltonian cycle in $(G^+)^2\/$ can be constructed (using $S$) which 
starts with $yx_1,x_1x_3\/$, ends with $yx_2\/$ and traverses $wx_4\/$ even if $w=z_2\/$ (see Observation (*)(i)-(ii) for $x_1x_3$). 

\vspace{3mm} {\bf (2.2)} \
Next assume that  $\{u_4, v_4\} = \{x_1, x_2\}\/$.

\vspace{1mm} 
Note that, in this case,  $d_G(x_2) > 2\/$ since $d_G(x_1) > 2\/$ can be assumed and $x_3, x_4\/$ are $2\/$-valent (note that the lemma is trivially true if 
$G$ is a $4$-cycle).

\vspace{1mm}
Consider the graph $G' = G -\{x_3, x_4\}\/$.

\vspace{1mm}
(a) Suppose $G'\/$ is $2\/$-connected. We shall apply Theorem \ref{2edge} to $G'\/$ with $x_1, x_2\/$ in place of $v, w\/$.

\vspace{1mm} 
(i) Suppose $G'\/$ has an $EPS\/$-graph $S' = E' \cup P'\/$ with $d_{P'}(x_i) =0\/$ for $i=1,2\/$. Let $E = E'  \cup \{yx_1x_4x_2y\}\/$ and 
$P = P' \cup \{x_1x_3\}\/$; this yields an $EPS\/$-graph $S = E \cup P\/$ of $G^+\/$ with  $d_P(x_1) =1\/$, $d_P(x_2) =0\/$, $d_P(x_3)=1\/$ and  
$d_P(x_4)=0\/$. Hence we may construct a hamiltonian cycle in $(G^+)^2\/$ containing the edges $x_1x_3\/$ and  $x_2x_4\/$ apart from $yx_1, yx_2\/$.

\vspace{1mm}
(ii) Suppose $G'\/$ has a $JEPS\/$-graph $S' = J' \cup E' \cup P'\/$ with $x_1, x_2\/$ being the only odd vertices of $J'\/$ and 
$d_{P'}(x_1) = 0 = d_{P'}(x_2)\/$. Let  $E = E' \cup (J' \cup \{x_1yx_2\}) \/$ and $P = P' \cup \{ x_1x_3, x_2x_4\}\/$. Then $S= E \cup P\/$ is 
an $EPS\/$-graph of $G^+\/$ with $d_P(x_1) = d_P(x_2) = d_P(x_3) = d_P(x_4) = 1\/$. Hence a hamiltonian cycle in $(G^+)^2\/$ containing the edges 
$yx_1, yx_2, x_1x_3\/$ and  $x_2x_4\/$ can be constructed.

\vspace{1mm}
(b) Finally assume that $G'\/$ is not $2\/$-connected. Then $G'\/$ is a non-trivial block chain. By Lemma \ref{flelemma}(ii) with $x_1=v\/$ and $x_2=w\/$, 
$G'\/$ has a $JEPS\/$-graph $S' = J' \cup E' \cup P'\/$  with $d_{P'}(x_1) =0= d_{P'} (x_2)\/$. As before, take $E = E'\cup (J' \cup \{x_1yx_2\})\/$ and 
$P = P' \cup \{x_1x_3, x_2x_4\}\/$. Then $S= E \cup P\/$ is an $EPS\/$-graph of $G^+\/$ with $d_P(x_1) = d_P(x_2) =  d_P(x_3)  = d_P(x_4) =1 \/$. Hence 
a hamiltonian cycle in $(G^+)^2\/$ containing the edges $yx_1, yx_2, x_1x_3\/$ and  $x_2x_4\/$ can be constructed.

\vspace{1mm} 
This completes the proof of the lemma.   \qed

\vspace{5mm} 
\noindent
{\bf Proof of Theorem \ref{dt}.}

\vspace{1mm} 
Let $G\/$ be a $2\/$-connected $DT\/$-graph and $A=\{x_1, x_2, x_3, x_4 \}\/$ be a set of four distinct vertices in $G\/$. It is easy to see that 
the theorem holds if $G\/$ is a cycle. Hence we also apply induction, apart from direct construction at the given graph. However, in general let $G^+\/$ 
be  defined as before.

\vspace{3mm} 
{\bf Case (A):}  $N(x_i) \subseteq V_2(G)\/$, $i=1,2,3,4\/$.

\vspace{1mm}
There exists a cycle $K^+\/$ in $G^+\/$ containing the vertices $y, x_1, x_2,  x_4\/$ (and possibly  $x_3\/$), assuming that $K^{+}$ is at least as long as 
any cycle containing $y, x_1, x_2,  x_3\/$. Assume $K^+\/$ is $W\/$-sound for $W=\{ y, x_1, x_2, x_3, x_4\}\/$. By Theorem \ref{thm4fh}, there exists 
a $W\/$-$EPS\/$-graph $S = E \cup P\/$ in $G^+\/$ with $K^+ \subseteq E\/$ (that is, $d_{P}(w) \leq 1\/$ for every vertex $w\/$ in $W\/$). Moreover 
$d_{P}(y) =0\/$ (since $y\/$ is  $2\/$-valent in $G^+\/$ and $K^+ \subseteq E\/$).

\vspace{1mm}
Since $N(x_i) \subseteq V_2(G)\/$ for $i=1, 2, 3, 4\/$, a hamiltonian cycle $C\/$ in $(G^+)^2\/$ can be constructed, and $C\/$ will contain $yx_1, yx_2\/$ 
and at least one edge of $G\/$ incident to $x_j\/$ for $j =3,4\/$. That is, $G^2\/$ contains a hamiltonian path as required (see Observation (*)(i)).

\vspace{3.88mm} 
{\bf Case (B):}  $N(x_i) \subseteq V_2(G)\/$, $i=1,2,3\/$ and  $N(x_4) \not \subseteq V_2(G)\/$; i.e.,  $d_G(x_4)=2\/$.

\vspace{3mm} 
Let $K^+\/$ be a cycle in $G^+\/$ containing $y, x_1, x_2, x_4\/$ and possibly $x_3\/$.

\vspace{2mm} 
{\bf (B)(1)} Suppose $x_3\/$ is not in $K^+\/$ (so, no cycle of $G^+\/$ contains $y\/$ and  $x_i\/$, $i=1, 2, 3, 4\/$).

\vspace{2mm} 
{\bf (a)} Suppose $\{u_4, v_4\} \neq \{x_1, x_2\}\/$.

\vspace{2mm} 
Then we may assume that $u_4  \not \in \{x_1,  x_2\}\/$.  Let $G' = G^+ - x_4\/$  and let $K' \subseteq G'\/$ be a cycle containing $y, x_1, x_2, x_3\/$.

\vspace{1mm} \hspace{5mm} {\bf (a1)}  
Suppose $G'\/$ is $2\/$-connected.

\vspace{1mm} 
Set $W' = \{y, x_1, x_2, x_3, u_4\}\/$ and suppose without loss of generality that $K'\/$ is $W'\/$-sound (i.e., $u_4 \in V(K')\/$ if $G'\/$ has a cycle 
containing all of $W'\/$). By Theorem \ref{thm4fh}, $G'\/$ has a $W'\/$-$EPS\/$-graph $S' = E' \cup P'\/$  with $K' \subseteq E'\/$ such that 
$d_{P'}(w) \leq 1\/$ for all $w \in W'\/$ with $d_{P'}(y)=0\/$. Take $E = E'\/$ and $P = P' \cup \{u_4x_4\}\/$. Then $S = E \cup P\/$ is an $EPS\/$-graph 
of $G^+\/$ with $K' \subseteq E\/$,  $d_P(y)=0, d_P(x_4)=1\/$ and  $d_P(w) \leq 1\/$ for $w \in W' -\{u_4\}\/$; and $d_P(u_4) \leq 2\/$. A careful 
examination of this case and Observation (*)(i)-(ii) show that a required hamiltonian cycle in $(G^+)^2\/$ can be constructed (note that $u_{4}x_{4}$ is 
a pendant edge of $S$).

\vspace{1mm} 
\hspace{5mm} {\bf (a2)} Suppose $G'\/$ is not $2\/$-connected.

\vspace{1mm} 
Then $G'\/$ is a non-trivial block chain. Let $B_y\/$ denote the block in $G'\/$ containing $K'\/$. Note that $u_4, v_4\/$ belong to different endblocks 
of $G'\/$. Let $z_4 \in \{u_4, v_4\}\/$ be a vertex in an endblock $B_1\/$ of $G'\/$ where $B_1 \neq B_y\/$. Further let $\widehat{G}\/$ denote the maximal 
block chain in $G'\/$ containing $B_1\/$ but no edges of $B_y\/$. Let $c_0 \in V(B_y) \cap V(\widehat{G})\/$ be a cutvertex of $G'\/$ (which is not 
a cutvertex of $\widehat{G}\/$).

\vspace{1mm}
Now replace $\widehat{G}\/$ in $G^+\/$ with a path $P_2 = z_4zc_0\/$ of length $2\/$ joining $z_4\/$ and $c_0\/$ and call the resulting graph $G^*\/$; 
$ z \not \in V(G)\/$. In so doing the cycle $K^+\/$ is transformed into the cycle $K^*\/$ in $G^*\/$ containing $P_2 \cup \{y, x_1, x_2, x_4\}\/$. Observe 
that $x_3 \not \in V(K^*)\/$; otherwise $K^*\/$ could be extended to become a cycle in $G^+\/$ containing $y, x_1, \ldots, x_4\/$ contrary to the supposition 
of this case.  Set $W = \{y, x_1, x_2, x_3, x_4\}\/$.  $K^*\/$ is $W\/$-sound in $G^*\/$;  by Theorem \ref{thm4fh}, $G^*\/$ contains a $W\/$-$EPS\/$-graph 
$S^* = E^* \cup P^*\/$ with $K^* \subseteq E^*\/$,  $d_{P^*}(y) = 0 = d_{P^*}(x_4) = d_{P^*}(z_4)\/$ and $d_{P^*}(w) \leq 1\/$ for all  $w\in W-\{y,x_4\}$.

\vspace{1mm}
Let $H = \widehat{G} \cup P_2\/$. Then $H\/$ is a $2\/$-connected graph and hence has a $[c_0; z_4]\/$-$EPS\/$-graph $S_H = E_H \cup P_H\/$ with 
$K_H \subseteq E_H\/$ where $K_H = (K^+ \cap \widehat{G}) \cup P_2\/$ (see Theorem~\ref{thm1f}).

\vspace{1mm} 
By taking $E = ((E^* \cup E_H) - (K^* \cup K_H)) \cup K^+\/$ and $P = P^* \cup P_H\/$ we have $S = E \cup P\/$ being  a $W\/$-$EPS\/$-graph of $G^+\/$   
with $K^+ \subseteq E\/$, $d_P(y)=0 = d_P(x_4)\/$, $d_P(w) \leq 1\/$ for all vertices $w\in W-\{y,x_4\}$, $d_P(z_4) \leq 1\/$ and $d_P(y_4) \leq 2\/$ where 
$y_4 \in N(x_4) -  z_4\/$. Hence a required hamiltonian cycle $H\/$ in $(G^+)^2\/$ can be constructed (as $\{u_4x_4, x_4v_4\} \subseteq K^+\/$); 
in particular $z_4x_4 \in E(H)\/$ (see Observation (*)(i)).

\vspace{2mm} 
{\bf (b)} Suppose $\{u_4, v_4\} = \{x_1, x_2\}\/$.

\vspace{1mm}   
Let $G' = G^+ - x_4\/$ (which is $2\/$-connected since $G\/$ is $2\/$-connected) and let  $K'\/$ be a cycle in $G'\/$ containing $y, x_1, x_2, x_3\/$.   
By Theorem \ref{thm2fh}, there exists an $[x_1; x_2, x_3]\/$-$EPS\/$-graph $S' = E' \cup P'\/$ in $G'\/$ with  $K' \subseteq E'\/$, $d_{P'}(w) \leq 1\/$ for 
$w \in \{ x_2, x_3\}\/$ and $d_{P'}(x_1) =0\/$. Let $E = E'\/$ and $P = P' \cup \{x_1x_4 \}\/$. Then we have an $EPS\/$-graph $S = E \cup P\/$ of $G^+\/$ 
with $K' \subseteq E\/$ and $d_P(x_i) \leq 1\/$ for $i=1,2,3,\/$ and $x_4\/$ is a pendant vertex in $S\/$.   Hence we can can construct a hamiltonian cycle 
in $(G^+)^2\/$ containing the edges  $yx_1, yx_2, x_1x_4\/$ and $x_3t_3\/$ where $t_3 \in N(x_3)\/$ since $N(x_3) \subseteq V_2(G)\/$ 
(see Observation (*)(i)-(ii)).

\vspace{3mm} 
{\bf (B)(2)} Suppose also $x_3\/$ is in $K^+\/$.

\vspace{1mm} 
Assume without loss of generality  that $K^+ = yx_1z_1 \cdots z_3x_3w_3 \cdots u_4x_4v_4 \cdots z_2x_2y\/$.

\vspace{2mm} 
{\bf (a)}   Suppose $u_4\neq x_3\/$.

\vspace{2mm}
Set $W = \{y, x_1, x_2, x_3, u_4\}\/$. Then $W \subseteq K^+\/$ and hence $K^+\/$ is $W\/$-sound. By Theorem \ref{thm4fh}, there is a $W\/$-$EPS\/$-graph 
$S = E\cup P\/$ in $G^+\/$ such that $K^+ \subseteq E\/$ and $d_{P}(w) \leq 1\/$ for every  $w \in W\/$. Then  it is possible to construct in $(G^+)^2\/$ 
a hamiltonian cycle  $C\/$ containing the edges $x_3w_3\/$ and $u_4x_4\/$ (recall that $x_4,  w_3\/$ are $2\/$-valent vertices in $G\/$) 
(see Observation~(*)(i)).

\vspace{2mm} 
{\bf (b)}  Suppose $u_4 = x_3\/$.

\vspace{2mm} 
(i) Suppose $v_4 \neq x_2\/$. We apply Theorem \ref{thm3fh} to $G^+\/$ to obtain  an $[x_3; x_1, x_2, v_4]\/$-$EPS\/$-graph $S= E \cup P\/$ with 
$K^+ \subseteq E\/$ and $d_{P}(x_3)=0\/$, $d_{P}(x_i)\leq 1\/$ for $i=1,2\/$, and $d_P(v_4) \leq 1\/$. Since $K^+ \subseteq E\/$ and $x_4 \in K^+\/$, 
we have  $d_{P}(x_4) =0\/$. We can construct a hamiltonian cycle $C\/$ in $(G^+)^2\/$ whose two edges incident to $x_i\/$ are edges of $G\/$ for $i=3\/$ or 
$i=4\/$, one of which is (without loss of generality) $x_3x_4\/$ (see Observation (*)(i)).

\vspace{2mm} 
(ii) Suppose $v_4 = x_2\/$. We operate analogously as in case (i) with an $[x_3; x_1, x_2, y]\/$-$EPS\/$-graph $S\/$ provided $x_3 \not \in N(x_1)\/$. 
However $S^2\/$ does not yield a hamiltonian cycle as required if $x_3 \in N(x_1)\/$. That is, $d_{G}(x_{3})=2$; $d_{G}(x_{4})=2$, and 
$N(x_{1})\subseteq V_{2}(G)$ by the assumptions. This is a special case of Lemma \ref{lemma4cd}. This finishes the proof of {\bf Case (B)}.  

\vspace{5mm}
{\bf Case (C):}  $N(x_i) \subseteq V_2(G)\/$, $i=1,2\/$ and $d_G(x_3) = 2 = d_G(x_4)\/$.

\vspace{3mm}
The proof of this case follows from  Lemma \ref{lemma4cd}.

\vspace{3mm} 
{\bf Case  (D):}  $N(x_1) \subseteq V_2(G)\/$ and  $N(x_2) \not  \subseteq V_2(G)\/$;  $d_G(x_2)=2\/$ follows.

\vspace{2mm} 
{\bf (D)(1) \ $N(x_4) \subseteq V_2(G)\/$. }

\vspace{1mm} 
There is a cycle $K^+$ in $G^+$ containing $y, x_1, x_2, x_3$ and also $x_4$ if a such a cycle exists. Recall that $x_3^*=x_3$ if $d_G(x_3) >2$ and 
$x_3^* = u_3=z_{3}\/$ if $d_G(x_3) =2\/$, and $N(x_2) = \{ u_2, v_2\}\/$ and assume that $v_2\/$ is in $K^+\/$. Let $x_3^-, x_3^+\/$ denote the predecessor, 
successor respectively, of $x_3\/$ in $K^+\/$, where we start the traversal of $K^{+}$ with the edge $yx_{1}$. We also note that 
$x_{3}^{*}=u_3=x_{3}^{-}$ and $v_3=x_{3}^{+}$ if $x_{3}\in V_{2}(G)$.

\vspace{3mm} 
{\bf (1.1)} Assume that $v_2 \not \in \{x_3, x_4\}\/$.

\vspace{2mm}  
{\bf (a)} $N(x_3)  \subseteq V_2(G)\/$.

\vspace{2mm} 
Let $W = \{y, x_1, v_2, x_3, x_4 \}\/$. Without loss of generality let $K^+\/$ be chosen such that it is $W\/$-sound, since 
$\{y, x_1, v_2, x_3\} \subseteq  K^+\/$  anyway, and possibly $x_4 \in K^+\/$. Let $S= E \cup P\/$ be a $W$-$EPS\/$-graph of $G^+\/$ with $K^+ \subseteq E\/$ 
(by Theorem \ref{thm4fh}). Observe that if $x_4 \not \in E\/$, then it is a pendant vertex in $S\/$; also $d_P(x_2) \leq 1\/$ automatically since 
$N(x_2) \not \subseteq V_2(G)\/$ and $x_2 \in K^+\/$. Now it is easy to construct a required hamiltonian cycle $C\/$ in $S^2\/$ having the required 
properties; we may assume that $x_3x_3^+ \in E(C)\/$ and $x_4w_4 \in E(G) \cap E(C)\/$, since $d_P(x_4) \leq 1\/$ and $N(x_4)\subseteq V_2(G)$. This is even 
true if $x_1=x_3^{-} \/$ since both $x_3\/$ and $x_3^+\/$ are $2\/$-valent in $G$ in this case (see Observation (*)(i)-(ii)).

\vspace{2mm}  
{\bf (b)} $N(x_3) \not \subseteq V_2(G)\/$; $d_G(x_3) =2\/$ follows.

\vspace{2mm} 
Set $W^+ = \{y, x_1, v_2, x_3^+, x_4 \}\/$.

\vspace{2mm} 
\hspace{5mm} {\bf (b1)} Suppose $x_4 = x_3^+\/$. Since $K^+ \supset \{y, x_1, v_2, x_3, x_3^+, x_4\}\/$, by Theorem \ref{thm3fh}, $G^+\/$ contains 
an $[x_4; y, x_1, v_2]\/$-$EPS\/$-graph $S= E \cup P\/$ with $K^+ \subseteq E\/$ such that $d_P(x_4)=d_P(x_3)=0$, $d_P(v_2) \leq 1\/$, $d_P(x_1) \leq 1\/$,
but also $d_P(x_2)\leq 1$. We obtain a hamiltonian cycle $C\subset S^2$ as required with $x_3x_3^+, x_4w_4\in E(G)\cap E(C)$ ($w_4\notin V(K^+)$ may hold, 
if $d_S(x_4)>2$). This covers also the case $x_1x_3\in E(K^+)$.

\vspace{2mm} 
\hspace{5mm} {\bf (b2)} Suppose $x_4 \neq x_3^+\/$.

\vspace{2mm}
\hspace{12mm} {\bf (b2.1)} Now assume that $K^+\/$ is $W^+\/$-sound.

\vspace{1mm} 
(i) Suppose $x_3^+\neq v_2$. Let $S = E\cup P$ be a $W^+\/$-$EPS\/$-graph of $G^+\/$ with $K^+ \subseteq E$, by Theorem \ref{thm4fh}. Then $S^2\/$ contains  
a hamiltonian cycle $C\/$ of $(G^+)^2\/$ as required, even if $x_1x_3\in E(K^+)$ and $d_G(x_3^+) >2\/$. In any case, also here $C\/$ can be constructed 
from $S\/$ such that $x_3x_3^+, x_4w_4 \in E(G) \cap E(C)\/$.

\vspace{1mm}
(ii) Suppose $x_3^+= v_2$. Hence $x_4\in K^+$, otherwise $|V(K^+)\cap W^+|=3$ and $K^+$ is not $W^+$-sound, a contradiction. Then $G^+\/$ contains 
an $[x_2; x_1, x_3^+,x_4]$-$EPS$-graph $S= E\cup P$ with $K^+ \subseteq E$, by Theorem \ref{thm3fh}. Hence we obtain a hamiltonian cycle $C\subset S^2$ 
as required with $x_3x_3^+, x_4w_4\in E(G)\cap E(C)$.
 
\vspace{2mm}
\hspace{12mm} {\bf (b2.2)} Assume that $K^+\/$ is not $W^+\/$-sound.

\vspace{1mm} 
(i) Suppose $|V(K^+)\cap W^+|>3$. Then there exists a cycle $K^*\/$ in $G^+\/$ containing $y, x_1, v_2, x_3^+, x_4\/$ but $x_3 \not \in K^*\/$; otherwise, 
we should have chosen $K^+ = K^*\/$ which is $W^+$-sound, a contradiction.

\vspace{1mm} 
First suppose $x_2v_2\in E(K^*)$. By Theorem \ref{thm3fh} (if $x_3^+\neq v_2$), Theorem \ref{thm2fh} (if $x_3^+=v_2$), there is 
an $[x_3^+; x_1, v_2, x_4]\/$-$EPS\/$-graph, $[x_3^+; x_1, x_4]\/$-$EPS\/$-graph, respectively, $S = E \cup P\/$ of $G^+\/$  with $K^* \subseteq E\/$. 
Note that either $x_3\/$ is a vertex in $E\/$, or else it is a pendant vertex in $S\/$. Also take note that $d_P(x_2) \leq 1\/$ and $d_P(v_2) \leq 1\/$. 
By Observation (*) (i)-(ii), $S^2\/$ has a hamiltonian cycle with the required properties.

\vspace{1mm} 
If $x_2v_2\notin E(K^*)$, then $x_2u_2\in E(K^*)$ and we proceed analogously as before using $u_2$ instead of $v_2$. Note that $u_2=x_4$ is not an obstacle 
(we use Theorem \ref{thm2fh}) because of $d_S(x_2)=2$ since $d_G(x_2)=2$ and $v_2$ is also in $K^*$ (thus $x_2v_2\notin E(S)$); and $v_2=x_3^+$ is not 
possible in this case.

\vspace{1mm} 
(ii) Suppose $|V(K^+)\cap W^+|=3$.

\vspace{1mm}
Hence $x_4\notin K^+$ and $v_2=x_3^+$. If $x_1x_3\notin E(K^+)$, then we set $W^*=\{y,x_1,x_3^-,x_3^+,x_4\}$. If $x_1x_3\in E(K^+)$ and $d_G(v_2)=2$, 
then we set $W^*=\{y,x_1,x_2,x_3^+,x_4\}$. In both cases $K^+$ is $W^*$-sound. By Theorem \ref{thm4fh}, $G^+$ contains a $W^*$-$EPS$-graph 
$S = E\cup P$ with $K^+ \subseteq E$. Observe that if $x_4 \not \in E\/$, then it is a pendant vertex in $S\/$. Now it is easy to construct a required 
hamiltonian cycle $C\/$ in $S^2\/$ having the required properties (see Observation (*)(i)-(ii)).

\vspace{1mm}
If $x_1x_3\in E(K^+)$ and $d_G(v_2)>2$, then we consider $G-x_3$. 

\vspace{1mm}
If $G-x_3$ is 2-connected, then we apply induction and get an ${\cal F}_4$ $x_1x_2$-hamiltonian path $P_1$ in $(G-x_3)^2$ containing edges 
$x_4w_4,x_3^+w_3^+\in E(G)$. Then $$P=P_1\cup \{w_3^+x_3,x_3x_3^+\}-\{x_3^+w_3^+\}$$ defines a hamiltonian path in $G^2$ as required. 

\vspace{1mm}
If $G-x_3$ is not 2-connected, then $x_1$ belongs to one endblock and $x_3^+,x_2$ to the other endblock of a non-trivial block chain $G-x_3$ because 
of the degree condition of $x_3, x_2$. Moreover $x_1$, $x_3^+$, and $x_2$ are not cutvertices of $G-x_3$. Depending on the position of $x_4$ in $G-x_3$ 
we construct a hamiltonian path $P$ in $G^2$ as in the preceding case applying either induction, or Theorem \ref{fletheorem4}, proceeding block after block. 
Since this procedure is straightforward we do not work out the details.

\vspace{3mm} {\bf (1.2)} $v_2  \in \{x_3, x_4\}\/$.

\vspace{3mm} \hspace{8mm}
{\bf (1.2.1)}  Assume that $v_2 = x_3\/$.

\vspace{1mm}
$G^+ - x_2u_2\/$ is a trivial or non-trivial block chain. Let $B_y\/$ denote the endblock (in $G^+ - x_2u_2\/$) containing the cycle $K^+\/$ and let 
$\widehat{G} = (G^+ -x_2u_2) - B_y\/$. 

Suppose first that $\widehat{G}\neq \emptyset$. Hence $\widehat{G}\/$ is a block chain in $G^+ -x_2u_2\/$ containing $u_2\/$ and $N(u_2)- x_2\/$. Let 
$t\/$ be the cutvertex of $G^+ - x_2u_2\/$ belonging to  $B_y\/$. That is, $\widehat{G} \cap B_y = \{t\}\/$.

\vspace{3mm} \hspace{15mm}
{\bf  (a)} Suppose $x_4 \not \in \widehat{G} -t \/$.

\vspace{1mm} 
Then $x_4\/$ is in $B_y\/$. Observe that $G-x_2\/$ is a block chain with $\widehat{G} \/$ being an induced subgraph of $G -x_2\/$ (note that $d_G(x_2) =2\/$).  
Since $B_y\/$ is $2\/$-connected, it contains a path $P(x_3, x_1)\/$ through $x_4\/$. It follows that $x_2, y \not \in P(x_3, x_1)\/$. Thus 
$P(x_3, x_1) \subseteq G-x_2\/$ with $x_4 \in P(x_3, x_1)\/$. Now, $P(x_2, x_1) = x_2x_3 P(x_3, x_1)\/$ is a path in $G - \widehat{G}\/$. Thus we may assume
that  $K^+ = yx_2P(x_2, x_1)x_1y \subseteq   B_y\/$ and thus  passes through $y, x_1, x_2, x_3^*, x_4\/$. By Theorem \ref{thm2fh}, $B_y\/$ has 
an $[x_3^*; x_1,  x_4]\/$-$EPS\/$-graph  $S_y = E_y \cup P_y\/$ with $K^+ \subseteq E_y\/$ and $d_{P_y}(x_2) =0\/$ (note that $d_{B_y}(x_2) =2\/$). Let 
the same $S_y\/$ denote an $[x_3^*; x_1]$-$EPS\/$-graph of $B_y\/$ if $x_4= x_3^*\/$ (i.e., $x_3x_4 \in E(K^+)\/$) (see Theorem~\ref{thm1f}).

\vspace{1mm} 
Since $x_4 \not \in \widehat{G}-t\/$, by Lemma \ref{flelemma}(i) or Theorem~\ref{thm1f} if $\widehat{G}$ is 2-connected, $\widehat{G}$ contains
an $EPS\/$-graph $\widehat{S} = \widehat{E} \cup \widehat{P}\/$ with $d_{\widehat{P}}(t) =0\/$ and $d_{\widehat{P}}(u_2) \leq 1\/$, provided 
$d_{\widehat{G}}(t)>1$. If $d_{\widehat{G}}(t)=1$, then either $\widehat{S}=\emptyset$ if $\widehat{G}=u_2t$ or by Lemma \ref{flelemma}(i), $\widehat{G}-t$ 
contains an $EPS\/$-graph $\widehat{S} = \widehat{E} \cup \widehat{P}\/$ with $d_{\widehat{P}}(t_1)\leq 1\/$ and $d_{\widehat{P}}(u_2) \leq 1\/$, where 
$t_{1}\in N_{\widehat{G}}(t)$.

\vspace{1mm}
Since $P_y \cap \widehat{P} = \emptyset\/$, by setting  $E = E_y \cup \widehat{E}\/ $ and $P = P_y \cup \widehat{P} \cup \{u_2x_2\}\/$, we obtain 
an $EPS\/$-graph $S = E \cup P\/$ of $G^+\/$ with $K^+ \subseteq E\/$, $d_P(x_2)=1\/$, $ d_P(x_3^*)=0\/$, $d_P(x_i) \leq 1\/$ for  $i=1, 4\/$ and 
$d_P(u_2) \leq 2\/$. 

If $|V(K^+)| \geq 6$, a required hamiltonian cycle in $S^2\/$ can be constructed (note that the cases $x_4 =t\/$ and $x_4 \neq t\/$ are treated 
simultaneously) (see Observation (*)(i)).   

\vspace{1mm} 
If, however,  $|V(K^+)| < 6$, i.e.,  $|V(K^+)| =5$, then $d_G(x_1)=2\/$ since $N(x_4) \subseteq V_2(G)\/$ and    $N(x_1) \subseteq V_2(G)\/$, which in turn 
implies $d_G(x_4)= d_G(x_3) =2\/$. Hence $G - \{x_3, x_4\}\/$ is a block chain $G^-\/$ with $x_1, x_2\/$ being pendant vertices of $G^-\/$. It follows that 
$(G^-)^2\/$ has a hamiltonian path $HP^-\/$  starting with $x_1w_1 \in E(G)\/$  and ending with $x_2u_2 \in E(G)\/$. Clearly, 
$$HP = (HP^- - \{x_1w_1, x_2u_2\}) \cup \{x_1x_4w_1, x_2x_3u_2\}$$ defines a required hamiltonian path in $G^2\/$.

\vspace{3mm} \hspace{15mm}
{\bf  (b)} Suppose $x_4  \in \widehat{G} -t\/$.

\vspace{1mm} 
In this case, we note that in $B_y$, the cycle $K^+$ can be assumed to traverse $y,x_1,t,x_3,x_2$ in this order; it also contains $x_3^*$ if $x_3$ is 
$2$-valent. As for $t \in V(K^+)\/$, see the preceding observation at the beginning of (a), with $t$ assuming the role of $x_4$. 

\vspace{1mm}  
Suppose $t\neq x_1$. By Theorem~\ref{thm2fh}, $B_y$ has an $[x_3^*;x_1,t]$-$EPS$-graph $S_y=E_y\cup P_y$ with $K^+ \subseteq E_y\/$ if $x_3^*\neq t$; by 
Theorem~\ref{thm1f}, $B_y$ has an $[x_3^*;x_1]$-$EPS$-graph $S_y=E_y\cup P_y$ with $K^+\subseteq E_y$ if $x_3^*= t$; and $d_{P_y}(x_2)=0$ since $d_G(x_2)=2$.

\vspace{1mm}
Suppose $t=x_1$. We let $S_y=E_y\cup P_y$ be an $[x_1; x_3^*]$-$EPS$-graph in $B_y$ with $K^+ \subseteq E_y$ by Theorem~\ref{thm1f}. Note that we set 
$x_3^*=x_2$ if $x_1x_3\in E(K^+)$.

\vspace{2mm} 
\hspace{18mm} {\bf (b1)} Assume  that  $x_4\/$ is a cutvertex in $\widehat{G}\/$.

\vspace{1mm} 
(i) Consider the  case $x_4\/$  is not incident to a bridge of $\widehat{G}$. Let $\widehat{G_1}\/$ and $\widehat{G_2}\/$ be defined by 
$\widehat{G} = \widehat{G_1} \cup \widehat{G_2}\/$ with $t, x_4 \in V(\widehat{G_1})\/$,  $x_4, u_2 \in V(\widehat{G_2})\/$ and  
$\widehat{G_1} \cap \widehat{G_2}= \{x_4\}\/$.

\vspace{1mm}
By Lemma \ref{flelemma}(i) or Theorem~\ref{thm1f}, $\widehat{G_i} \/$ has an $EPS\/$-graph $\widehat{S_i} = \widehat{E_i} \cup \widehat{P_i}\/$ with 
$d_{\widehat{P_i}}(x_4) =0 \/$ for $i=1, 2\/$, $d_{\widehat{P_1}}(t) \leq 1\/$ and $d_{\widehat{P_2}}(u_2) \leq 1\/$.

\vspace{1mm}
By taking $E = E_y \cup \widehat{E_1} \cup \widehat{E_2}\/$, $P = P_y \cup \widehat{P_1} \cup \widehat{P_2}\/$, we have an $EPS\/$-graph $S = E \cup P\/$ of  
$G^+\/$ with $d_{P}(x_2) =0 = d_P(x_4)\/$,  and $d_{P}(x_1) \leq 1\/$ and $d_P(x_3^*) \leq   1\/$; $d_P(t) \leq 2\/$ by construction, provided $t \neq x_1\/$. 
Moreover, if $t=x_3^*$, we have $d_{P_y}(x_3^*)=0\/$ and hence $d_P(x_3^*) \leq 1\/$ because of $d_{\widehat{P_1}}(x_3^*) \leq 1\/$; and $d_P(x_1) \leq 1\/$. 
Also if $t = x_1\/$, we have $d_{P_y}(x_1)=0\/$ and hence $d_P(x_1) \leq 1\/$ because of $d_{\widehat{P_1}}(t)\leq 1\/$; and $d_P(x_3^*)\leq 1.$ Hence 
a required hamiltonian cycle in $S^2\/$ can be constructed (the various construction details are straightforward and are thus omitted).   

\vspace{2mm} 
(ii) Now suppose $x_4\/$ is incident to a bridge $f\/$ of $\widehat{G}$ and $|V(K^+)| > 4\/$. In this case, we delete  $f\/$  and thus split $\widehat{G}\/$ 
into two block chains $\widehat{G_1}\/$ and $\widehat{G_2}\/$ with $t \in \widehat{G_1}\/$, $u_2 \in \widehat{G_2}\/$ and $x_4\/$ is either in $\widehat{G_1}$
or in $\widehat{G_2}$. By Lemma \ref{flelemma}(i) or Theorem \ref{thm1f}, $\widehat{G_i}$ has an $EPS$-graph $\widehat{S_i}=\widehat{E_i}\cup \widehat{P_i}$ 
with $d_{\widehat{P_1}}(t)\leq 1$, $d_{\widehat{P_2}}(u_2)\leq 1$ and $d_{\widehat{P_i}}(x_4)\leq 1$ for some $i\in \{1, 2\}$. Note that $\widehat{S_i}=G_2$ 
if $G_i=K_2$; or $\widehat{S_i}=\emptyset$ if $G_i=t$ or $G_i=u_2$. Proceeding similarly to case (i) let $E = E_y \cup \widehat{E_1} \cup \widehat{E_2}$ and
$P = P_y \cup \widehat{P_1} \cup (\widehat{P_2} \cup \{u_2x_2\})\/$. Then we have an $EPS\/$-graph $S= E \cup P\/$ of $G^+$.

\vspace{1mm} 
Because of the choice of $S_y\/$ in the cases $t \notin \{x_1,x_3^*\}\/$, $t=x_3^*$, and $t = x_1$, we have in any case, $d_{P}(x_1) \leq 1\/$, 
$d_{P}(x_2) =1\/$, $d_P(x_3^*) \leq 1\/$ and $x_4\/$ is  either a pendant vertex in $S\/$ or $d_P(x_4) =0\/$ (which occurs when $d_G(x_4) >2\/$). 

\vspace{1mm}
By a similar argument as in case (i), we conclude that in all cases $S^2\/$ contains a hamiltonian cycle with the required properties unless 
$d_{P}(x_2) =1$, $d_P(x_3^*)=1$ and $x_3^*=x_3=t$. In this case there exists a cycle containing $y,x_1,x_2,x_3,x_4$, a contradiction to the choice of $K^+$.

\vspace{2mm} 
(iii) Now suppose $x_4\/$ is incident to a bridge $f\/$ of $\widehat{G}$ and $|V(K^+)| = 4\/$. It follows that $t=x_{1}$ and therefore $G' = G - x_3\/$ is 
a non-trivial block chain containing $f$ as a bridge. Hence $(G-x_3)^2\/$ contains a hamiltonian path $HP'\/$ starting with an edge $x_1w_1 \in E(G)\/$ 
and ending with $u_2x_2 \in E(G)\/$ and containing  $f\/$. It follows that $$HP = (HP' - x_1w_1) \cup \{x_1x_3w_1\}$$ yields a hamiltonian path in $G^2\/$ 
as required if $f\neq x_{1}x_{4}$. However, if $f=x_{1}x_{4}$, then we set $$HP = (HP' - x_2u_2) \cup \{x_2x_3u_2\}.$$

\vspace{2mm} 
\hspace{18mm} {\bf (b2)} Hence assume that $x_4\/$ is not a cutvertex in $\widehat{G}\/$.

\vspace{1mm}  
Suppose first that $x_4\/$ is contained in a $2\/$-connected block $B\/$ in $\widehat{G}\/$. Further, let  $c_1, c_2 \/$ be two vertices in $B\/$ which are 
also cutvertices of $\widehat{G}\/$ if $B\/$ is not an endblock of $\widehat{G}\/$. If, however, $B\/$ is an endblock of $\widehat{G}\/$, then let $c_1\/$ 
be the unique cutvertex of $\widehat{G}\/$ in $B\/$, and let $c_2 \in \{t, u_2\}\/$ depending on which of the endblocks of $\widehat{G}\/$ is $B\/$.         
If $x_{4}\neq c_{2}$, we apply Theorem \ref{thm2fh} to $B\/$ to obtain an $[x_4; c_1, c_2]\/$-$EPS\/$-graph of $B\/$; if $x_{4}=c_{2}$ (which means 
$x_{4}=u_{2}$), then we apply Theorem \ref{thm1f} to $B\/$ to obtain an $[x_4; c_1]\/$-$EPS\/$-graph of $B\/$. In both cases by using Lemma \ref{flelemma}(i), 
extend these $EPS\/$-graphs to an $EPS\/$-graph $\widehat{S} = \widehat{E} \cup \widehat{P}\/$ of $\widehat{G}\/$ with $d_{\widehat{P}}(t) \leq 1\/$,
$d_{\widehat{P}}(u_2) \leq 1\/$, and $d_{\widehat{P}}(x_4) =0\/$.

\vspace{1mm}
Setting $E = E_y \cup \widehat{E}\/$ and $P = P_y \cup \widehat{P}\/$, we obtain an $EPS\/$-graph $S= E \cup P\/$ of $G^+\/$ with 
$d_{P}(x_2) =0 = d_P(x_4)\/$, $d_P(x_3^*)\leq 1\/$ and $d_{P}(x_1) \leq 1\/$.

\vspace{2mm}
Hence assume that $x_4\/$ is not contained in a $2$-connected block. That is, $x_4\/$ is a pendant vertex  in $\widehat{G}\/$. In this case, $x_4 = u_2\/$. 
We apply Lemma \ref{flelemma}(i) to obtain an $EPS$-graph $\widehat{S} = \widehat{E} \cup \widehat{P}$ of $\widehat{G}\/$ with $d_{\widehat{P}}(t) \leq 1$, 
and $d_{\widehat{P}}(x_4) \leq 1\/$ if $\widehat{G}\neq x_4t$. If $\widehat{G}=x_4t$, then $\widehat{S}=\widehat{G}$. Setting $E = E_y \cup \widehat{E}$ 
and $P = P_y \cup \widehat{P}$, we obtain an $EPS$-graph $S= E \cup P\/$ of $G^+\/$ with $d_{P}(x_2) =0 \/$ and $d_P(x_3^*)\leq 1\/$, $d_{P}(x_1) \leq 1$, 
$d_P(x_4) = 1\/$ and $x_4\/$ is a pendant vertex in $S$.

\vspace{1mm}
In any of these cases, $S^2\/$ contains a hamiltonian cycle $C\/$ with the required properties (note that $x_3x_2 \in E(C)\/$ because 
$d_E(x_2) = d_{G^+}(x_2) -1 =2\/$); see Observation(*)(i)-(ii).

Finally if $\widehat{G}=\emptyset$, we find $S_{y}$ as in Case {\bf (1.2.1)(a)} and construct a hamiltonian cycle as required using $S_{y}$ only.

\vspace{3mm}  
\hspace{8mm} {\bf (1.2.2)}  Assume that $v_2 = x_4\/$.

\vspace{1mm} 
Recall that the cycle $K^+$ in $G^+$ contains $y,x_1,x_2,x_3,v_2$. Therefore $$K^+=yx_1\ldots z_3x_3w_3\ldots z_4x_4x_2y.$$ Consider the graph 
$G'=G^+ - x_2u_2$.

\vspace{2mm} 
\hspace{15mm} {\bf Case (a)} \ $G'\/$ is $2\/$-connected.

\vspace{1mm}  
Suppose $x_3, x_4\/$ are adjacent in $K^{+}\/$. Then apply Theorem \ref{thm1f} to obtain  an $[x_4; x_1]\/$-$EPS\/$-graph $S = E \cup P\/$ of $G'\/$ with 
$K^+ \subseteq E\/$. Suppose $x_3, x_4\/$ are not adjacent in $K^{+}\/$. Then apply  Theorem \ref{thm2fh} to obtain an $[x_1; x_3^*, x_4]\/$-$EPS\/$-graph 
$S = E \cup P\/$ of $G'\/$ with $K^+ \subseteq E\/$. In either case, a required hamiltonian cycle in $S^2\/$ can be constructed (setting $x_3^* = x_3^+=w_3$
if $x_1x_3 \in E(K^+)\/$).

\vspace{2mm} 
\hspace{15mm} {\bf Case (b)} \ $G'\/$ is not $2\/$-connected.

\vspace{1mm}
Then $G'\/$ is a non-trivial block chain. As before, let $B_y\/$ denote the endblock in $G'\/$ containing $y\/$ (and hence containing the cycle $K^+\/$). 
Set $\widehat{G} = G' -B_y\/$ which is a trivial or non-trivial block chain; $\widehat{G}\neq \emptyset$ in any case. It follows that 
$B_y \cap \widehat{G} = \{t\}\/$ and $t\/$ is a cutvertex of $G'\/$. By Theorem \ref{thm1f} or Lemma \ref{flelemma}(i), $\widehat{G}\/$ has an $EPS\/$-graph 
$\widehat{S}=\widehat{E}\cup \widehat{P}$ with $d_{\widehat{P}}(t)\leq 1$ if $\widehat{G}\neq u_2t$. If $\widehat{G}=u_2t$, then $\widehat{S}=\widehat{G}$.

\vspace{2mm}
(i) Suppose $t = x_4\/$.

\vspace{1mm} 
Then $G'' = G^+ - x_2x_4\/$ is $2\/$-connected. Replace in $K^+\/$ the edge $x_4x_2\/$ with a path $P(x_4, x_2)\/$ in $\widehat{G} \cup \{u_2x_2\}\/$ 
to obtain  the cycle $K''\/$.  Since $\{y, x_1, x_3, x_4, u_2, x_2\} \subseteq V(K'')\/$, we may apply Theorem \ref{thm3fh} to obtain 
an $[x_4; x_1, u_2, x_2]\/$-$EPS\/$-graph $S'' = E'' \cup P''  \subseteq G''\/$ if $x_3\/$ and $x_4\/$ are adjacent in $K''\/$, or to obtain 
an $[x_1; x_3^*, x_4, u_2]\/$-$EPS\/$-graph $S'' = E'' \cup P''\subseteq G''$ if $x_3\/$ and $x_4\/$ are not adjacent in $K''\/$ (setting $x_3^* = x_3^+=w_3$ 
if $x_1x_3 \in E(K'')\/$).  In both cases, $K'' \subseteq E''\/$. A required hamiltonian cycle in $(S'')^2$ can be constructed (since the situation is 
similar to Case \textbf{(a)} above); see Observation (*)(i).

\vspace{1mm}
(ii) Suppose  $t = x_3\/$.

\vspace{1mm} 
We apply Theorem \ref{thm2fh} to $B_y\/$ to obtain an $[x_3; x_1, x_4]\/$-$EPS\/$-graph $S_y = E_y \cup P_y\/$ of $G'$ with $K^+ \subseteq E_y\/$. Note that 
$N(x_3) \subseteq V_2(G)\/$ in this case.

\vspace{1mm}
(iii) Suppose $t =x_1\/$ and $x_1x_3x_4 \not \subseteq K^+\/$.

\vspace{1mm} 
We set $x_3^*=x_3^+=w_3$, if $x_1x_3 \in E(K^+)$. We apply Theorem \ref{thm2fh} to $B_y\/$ again to obtain an $[x_1; x_3^*, x_4]\/$-$EPS\/$-graph 
$S_y = E_y \cup P_y\/$ of $G'$ with $K^+ \subseteq E_y\/$.

\vspace{1mm}
In the cases (ii) and (iii), we let $E = \widehat{E} \cup E_y\/$, $P = \widehat{P} \cup P_y\/$ and obtain an $EPS\/$-graph $S= E \cup P\/$ of 
$G^+ -x_2u_2$ with $K^+ \subseteq E\/$, $d_P(x_2)=0 \/$ and $d_P(w) \leq 1\/$ for $w \in \{x_1, x_3^*, x_4\}\/$. Hence a required hamiltonian cycle 
in $S^2\/$ can be constructed; see Observation (*)(i).

\vspace{1mm}
(iv) Suppose $t \not \in \{x_1, x_3, x_4\}\/$.

\vspace{1mm} 
We set $x_3^*=x_3^+=w_3$ if $x_1x_3\in E(K^+)$. Note that $d_{B_y}(x_2)=2$. Let $W=\{y,x_1,x_3^*,x_4,t\}$; $W-\{t\}\subseteq V(K^+)$. 

\vspace{1mm}
First suppose that $|W|=5$. If $K^+$ is not $W$-sound, then there is a cycle $K'\subseteq B_y$ containing all vertices of $W$, in which case we apply 
Theorem \ref{thm3fh} to $B_y\/$ to obtain an $[x_3^*; x_1, x_4, t]\/$-$EPS\/$-graph $S_y = E_y \cup P_y\/$ with $K' \subseteq E_y\/$. Note that, if 
$x_3\notin K'$, then either $x_3$ is a pendant vertex in $S_y$ or $d_{P_y}(x_3)=0$ and $x_3^*x_3\in E(E_y)$. If $K^+\/$ is $W\/$-sound, then we set $K'=K^+$
and apply Theorem \ref{thm4fh} to $B_y\/$ to obtain a $W\/$-$EPS\/$-graph $S_y = E_y \cup P_y\/$ with $K' \subseteq E_y\/$. 

\vspace{1mm}
Suppose $|W|=4$. Hence $x_3^*\neq x_3$. 

\vspace{1mm}
If $x_3^*=x_4\neq t$, then $N_G(x_3)=\{x_1,x_4\}$ and $B_y-x_3$ is 2-connected: for, there are two internally disjoint paths from $t$ to $K^+$, and the 
endvertices of these paths in $K^+$ are $x_1$ and $x_4$ since $x_3,x_2\in V_2(G)$. Thus $B_y-x_3$ contains a cycle $K'$ containing $y,x_1,x_2,x_4,t$. Hence 
we apply Theorem \ref{thm3fh} to $B_y-x_3$ to obtain an $[x_1; x_2, x_4, t]$-$EPS$-graph $S'_y = E'_y \cup P'_y\/$ with $K' \subseteq E'_y\/$. Let $E_y=E'_y$ 
and $P_y=P'_y\cup \{x_1x_3\}$. Thus we have an $EPS$-graph $S_y = E_y \cup P_y\/$ of $B_y$ with $K' \subseteq E_y\/$. Moreover $d_{P_y}(x_1)=1$, 
$d_{P_y}(x_2)=0$, $d_{P_y}(x_3)=1$, $d_{P_y}(x_4)\leq 1$, $d_{P_y}(t)\leq 1$, and $x_3$ is a pendant vertex in $S_y$. 

\vspace{1mm}
If $x_3^*=t\notin\{x_1,x_3,x_4\}$, then we set $K'=K^+$ and apply Theorem \ref{thm2fh} to $B_y$ to obtain an $[x_3^*;x_1,x_4]$-$EPS$-graph $S_y=E_y\cup P_y$ 
with $K'\subseteq E_y$.

\vspace{1mm} 
In all cases, we let $E = \widehat{E} \cup E_y\/$, $P = \widehat{P} \cup P_y\/$ and obtain an $EPS\/$-graph $S= E \cup P\/$ of $G^+ -x_2u_2\/$ with 
$K' \subseteq E\/$, $d_P(x_2)=0 \/$ and $d_P(w) \leq 1\/$ for $w \in \{x_1, x_3^*, x_4\}\/$ (even if $x_1x_3 \in E(K^+)\/$ or $x_3^* = x_3^+ =t \neq x_3\/$).
Hence a required hamiltonian cycle in $S^2\/$ can be constructed; see Observation (*)(i)-(ii).

\vspace{1mm}
(v) Suppose $t =x_1\/$ and $x_1x_3x_4  \subseteq K^+\/$.

\vspace{1mm} 
Note that $K^+ = y x_1 x_3 x_4x_2y\/$.  Let $G_3 = B_y -\{y, x_2\}\/$; it is $2\/$-connected if $d_G(x_4) >2\/$, or else it is a path $x_1x_3x_4\/$. 
We have $G- x_2x_4  = G_3 \cup (\widehat{G} \cup \{u_2x_2\})\/$ with $t = G_3 \cap \widehat{G}\/$. Consequently, 
$$ \widetilde{G} := \widehat{G} \cup \{u_2x_2\} =  G - (\{x_2x_4\} \cup G_3).$$
By Corollary \ref{flecor}(ii), $\widetilde{G}^2\/$ has a hamiltonian path $P_{1,2}\/$ starting with $x_1w_1 \in E(G)\/$ and ending with $u_2x_2\/$.  
If $G_3\/$ is $2\/$-connected, then $G_3 - x_3\/$ is a (trivial or non-trivial) block chain and thus $(G_3 - x_3)^2\/$ has a hamiltonian path $P_{4,1}\/$ 
starting in $x_4=v_2$ and ending with an edge $s_1x_1\in E(G)$ (using Theorem \ref{fletheorem4}(ii) if $G_3-x_3$ is 2-connected, Corollary \ref{flecor}(ii) 
if $G_3-x_3$ is a non-trivial block chain, and $P_{4,1}=G_3-x_3$ if $G_3-x_3=x_1x_4$). Set $$P(x_1, x_2) = x_1x_3w_1 (P_{1,2} - \{x_1, x_2\}) u_2x_4x_2$$ 
if $G_3 = x_1x_3x_4\/$; and $$ P(x_1, x_2) =  x_1x_3x_4 (P_{4,1} - x_1)s_1w_1 (P_{1,2}-x_1)$$ if $G_3\/$ is $2$-connected. In both cases, $P(x_1, x_2)\/$ is 
a ${\cal F}_4\/$ $x_1x_2\/$-hamiltonian path in $G^2\/$. This finishes the proof of Case {\bf (D)(1)}.

\vspace{3mm}
Since the case $N(x_3) \subseteq V_2(G)\/$ is analogous to the Case \textbf{(D)(1)}, we are left with the following case.

\vspace{3mm}
{\bf (D)(2)} $N(x_i) \not \subseteq V_2(G)\/$ for $i=3, 4\/$.

\vspace{3mm}
However, the proof of this case follows from Lemma \ref{lemma4cd}. This finishes the proof of Case {\bf (D)}.

\vspace{5mm}
{\bf Case (E):}  $N(x_1) \not \subseteq V_2(G)\/$ and  $  N(x_2) \not \subseteq V_2(G)\/$.

\vspace{2mm}
Then $d_G(x_1) = 2 = d_G(x_2)\/$.

\vspace{2mm}  
Let $K^+\/$ be a cycle containing the vertices $y, u_1, u_2, x_3\/$ and possibly $x_4\/$ where we  assume that $$K^+= yx_1u_1 \cdots x_3 \cdots u_2x_2y.\/$$

\vspace{3mm} {\bf (E)(1) } Suppose $x_4\/$ is not in any cycle containing $y\/$ and $x_3\/$.

\vspace{3mm}
{\bf (1.1)} \ $d_G(x_3) > 2, \ d_G(x_4) > 2\/$.

\vspace{2mm}
{\bf (a)}  Suppose $x_3 \not \in \{ u_1, u_2\}\/$.

\vspace{2mm}
Set $W = \{ y, u_1, u_2, x_3, x_4\}$. By supposition, $K^+$ is $W$-sound. By Theorem \ref{thm4fh}, we have an $EPS$-graph $S= E\cup P$ with $K^+\subseteq E$
and $d_P(w) \leq 1\/$ for every $w \in W\/$. In this case  a required hamiltonian cycle $C\/$ in $S^2\/$ can be constructed (taking note that $x_1, x_2\/$ 
are $3\/$-valent in $G^+\/$, and that $x_ix_4 \in E(P)\/$, $i  \in \{1, 2\}\/$, does not constitute an obstruction in the construction of $C\/$).

\vspace{2mm}
{\bf (b)} Suppose $x_3 = u_1\/$.

\vspace{1mm} 
Note that if $x_4  \not \in N(x_1) \cup N(x_2)\/$, then we are back to case (a) with $x_3\/$ and $x_4\/$ changing roles. Hence we have $x_4\in \{v_1, v_2\}$. 
Also, $x_4 =v_1=v_2\/$ cannot hold; otherwise,  $d_G(x_4) > 2\/$ and $x_i \in N(x_4)\/$, $d_G(x_i) =2\/$, $i=1, 2\/$ imply the existence of an $x_4q_3\/$-path 
$P(x_4, q_3) \subset G\/$ with $q_3 \in V(K^+)\/$ and  $(P(x_4, q_3) - q_3 ) \cap K^+ =\emptyset\/$,  yielding in turn  a cycle containing $y, x_3, x_4\/$ 
contradicting \textbf{E(1)}. By the same token, $x_3=u_1=u_2\/$ cannot hold. 

\vspace{2mm} 
\hspace{8mm} {\bf (b1)} $x_4 = v_2\/.$

\vspace{2mm}
Consider $G^- = G - \{x_1u_1, x_2u_2\}\/$.

\vspace{1mm}
Note that $x_3, x_4\/$ belong to different components of  $G^-\/$; otherwise there is a path  $P_0\/$ in $G^-\/$ joining $x_3\/$ and $x_4\/$ implying that  
$C_0 = P_0 x_4 x_2 yx_1x_3\/$ is a cycle in $G^+\/$ with $y, x_3, x_4 \in V(C_0)\/$,  a contradiction to the supposition. Since $G\/$ is $2$-connected, $G^-$ 
contains precisely two components $G_3^-\neq K_1$ and $G_4^-$ containing $x_3, x_4$, respectively. Clearly $x_2\in V(G_4^-)$. We also have  $x_1\in V(G_4^-)$ 
because $P_0$ as above does not exist.

\vspace{1mm} 
Observe that $G_4^-$ and $G_3^-$ are (trivial or non-trivial) block chains in which $x_1, x_2 \in V(G_4^-)$ and $x_3, u_2 \in V(G_3^-)$ are not cutvertices.
Thus  $G^+-\{x_1u_1, x_2u_2\}\/$ is a disconnected graph with two components $G_3 = G_3^-\/$ (which contains $x_3=u_1\/$ and $u_2\/$)   and $G_4\/$ 
(which contains $y, x_1, v_1, x_4=v_2\/$ and $x_2\/$).

\vspace{1mm} 
Note that in $G_4\/$, there is a cycle $C^+$ containing $y, x_1, v_1, x_4, x_2$, implying that $G_4\/$ is $2\/$-connected, whereas $G_3\/$ is a block chain. 
By Theorem \ref{thm1f}, let $S_4=E_4\cup P_4$ be a $[v_1; x_4]$-$EPS$-graph in $G_4$ with $C^+\subseteq E_4$, $d_{P_4}(v_1)=0=d_{P_4}(x_1)=d_{P_4}(x_2)$ and
$d_{P_4}(x_4)\leq 1$. By Lemma \ref{flelemma}(i) or Theorem \ref{thm1f} (respectively depending on whether $G_3$ has a cutvertex or $G_3$ is $2$-connected), 
there is an $EPS\/$-graph $S_3 = E_3 \cup P_3\/$ in $G_3\/$ such that $d_{P_3}(x_3)=0\/$ and $d_{P_3}(u_2) \leq 1\/$. Taking $E= E_3 \cup E_4\/$ and
$P=P_3\cup\{x_1x_3\}\cup P_4$, we have an $EPS$-graph $S=E\cup P$ of $G^+\/$ with $C^+ \subseteq E\/$ and $d_P(v_1)=0= d_P(x_2)\/$, $d_P(x_1) =1 = d_P(x_3)$
and $d_P(x_4) \leq 1\/$.

\vspace{2mm} 
Note that in this case, since $d_G(x_3), \ d_G(x_4) > 2\/$,  $d_P(x_2)=0$, $d_P(x_3) =1\/$ and $d_P(x_4) \leq 1\/$, it is straightforward that one can obtain   
a required hamiltonian cycle of $(G^+)^2\/$.

\vspace{3mm} 
\hspace{8mm} {\bf (b2)}  $x_4 = v_1$.

\vspace{2mm}
Let $G' = G^+ -x_1x_3\/$, and we may assume that a cycle  $K' = yx_1x_4\cdots v_2 x_2y \subseteq G'$ exists. Note that in $G$ we have two 
internally disjoint paths $x_1x_3\cdots u_2x_2$ and $x_1x_4\cdots v_2x_2$. This is in line with the notation of $K^+$ above.

\vspace{2mm}  
\hspace{15mm} {\bf (b2.1)} Suppose $G'\/$ is $2\/$-connected.

\vspace{2mm} 
Take $W= \{y, x_3, x_4, v_2, x_2\}\/$. Then $K'\/$ is $W$-sound in $G'\/$  since $v_2 \neq x_4\/$ (see the observation in \textbf{(b)}). Let $S = E \cup P\/$ 
be a $W\/$-$EPS\/$-graph of $G'\/$ (and hence a $W\/$-$EPS\/$-graph of $G^+\/$) with $K' \subseteq E\/$ and $d_P(w) \leq 1\/$ for every $w \in W\/$. Since 
$d_G(x_3)>2$, $d_G(x_4)>2$, a hamiltonian cycle in $S^2\/$ can be constructed containing $x_1y, x_2y$ and  $x_iz_i$ where $z_i \in N_G(x_i)$, $i=3,4$.

\vspace{2mm} 
\hspace{15mm}
{\bf (b2.2)} Suppose $G' \/$ is not $2\/$-connected.

\vspace{2mm} 
By symmetry, $ G^+ -x_1x_4\/$ is also not $2\/$-connected. Then $G'\/$ is a block chain with endblocks  $B_3, B'\/$, with $x_3 \in B_3\/$ and $K' \subset B'$ 
and $x_1\/$ and $x_3\/$ are not cutvertices of $G'\/$. Furthermore, let $c\/$ denote the cutvertex of $G'\/$ which belongs to $B'\/$; $c\neq x_4$ (otherwise, 
$G^+\/$ contains a cycle through $y, x_3, x_4\/$).

\vspace{2mm}
Set $G_0 = G' -B'\/$. Note that  $x_3, c\/$ are vertices in $G_0\/$ and are not cutvertices of $G_0\/$. By Lemma \ref{flelemma}(i) or Theorem \ref{thm1f} 
(depending on whether $G_0$ has a cutvertex or not), $G_0$ contains an $EPS\/$-graph $S_0 = E_0 \cup P_0\/$ with $d_{P_0}(c) \leq 1\/$ and $d_{P_0}(x_3) =0\/$ 
($B_3\/$ is $2\/$-connected because $d_{B_3}(x_3) >1\/$).

\vspace{3mm}
(i) Suppose $c \not \in \{v_2, x_2\}\/$. Let $W' = \{y, x_4, c, v_2, x_2\}\/$. $B' \supseteq K' \supset (W' - c)\/$ in any case. So, $K'\/$ is $W'$-sound,  
or there is a cycle $K'' \supset W'\/$ with $B' \supseteq K''\/$, in which case $K''\/$ is $W'\/$-sound in $B'\/$.

\vspace{2mm} 
(ii) Now suppose $c=x_2$. Set $W' = \{y, x_1, x_2, v_2, x_4\}\/$ and observe that $K'\/$ is $W'$-sound in $B'\/$ again.

\vspace{2mm} 
In both cases, we obtain by Theorem \ref{thm4fh} an $EPS$-graph $S'=E'\cup P'$ of $B'$ with $K'\subseteq E'$ or $K''\subseteq E'$, and $d_{P'}(w) \leq 1\/$ 
for every $w\in W'$. Note that if $c\notin \{v_2,x_2\}$, $c\notin K'$ and $x_2v_2\notin E(K'')$, or if $c=x_2$, then $d_{P'}(x_2)=0$ because $d_{B'}(x_2)=2$.

\vspace{2mm}
Set $E=E_0\cup E'$, $P=P_0\cup P'$ to obtain an $EPS$-graph $S=E\cup P$ of $G^+$ with $K^*\subseteq E$ where $K^*\in\{K',K''\}$, $d_P(x_3)=0$, $d_P(z)\leq 1$ 
for every $z \in \{y, x_4, v_2, x_2\}\/$, and  $d_P(c) \leq 2\/$ if $c \not \in \{v_2, x_2\}\/$, and $d_P(c)\leq 1\/$ if $c=x_2\/$. Also, $d_P(x_1) =0\/$ 
since $x_1x_3 \not \in E(S)$. Since $N(x_i) \subseteq V_2(G)$, $i=3,4$ and $d_P(x_3)=0$, $d_P(x_4)\leq 1$, a hamiltonian cycle in $S^2$ containing the edges
incident to $y$ and containing edges $x_iz_i$, can be constructed, where $z_i\in N_G(x_i)$, $i=3,4$. Observe that $d_P(v_2)=d_P(x_2)=1$ does not create 
any obstacle.

\vspace{2mm}
(iii) Suppose $c=v_2\/$. In this case, by Theorem~\ref{thm2fh} we take in $B'\/$ a $[v_2; x_2, x_4]\/$-$EPS\/$-graph and proceed as in case (i).

\vspace{3.88mm}
{\bf (1.2)} $d_G(x_3)>2, d_G(x_4)=2$.

\vspace{1mm} 
Let $K'$ be a cycle in $G^+\/$ containing the vertices $y, x_1, w_1,  u_4, x_4, v_4, w_2, x_2\/$ in this order where $w_i \in \{u_i, v_i\}\/$, $i =1, 2\/$.

\vspace{2mm}
{\bf (a)}   $x_4 \notin \{w_1, w_2\}\/$

\vspace{2mm} 
\hspace{5mm} {\bf (a1)} Suppose $v_4 \not \in N(x_2)\/$. Note that in this case $|V(K')|>6\/$.

\vspace{1mm}  
Set $W=\{y,w_1, w_2, x_3,  v_4\}\/$ and observe that $|W|=5$ and $|K'\cap W|\geq 4$.

\vspace{1mm}  
Suppose $K'$ is $W$-sound in $G^+$. Then by Theorem \ref{thm4fh}, $G^+$ has a $W$-$EPS$-graph $S=E\cup P$ with $K'\subseteq E$ and $d_P(y)=0=d_P(x_4)$. 
Moreover, for $i=1,2$, we have $d_P(x_i)\leq 1$ since $d_G(x_i)=2$. Hence we can construct a hamiltonian cycle in $S^2$ having the required properties.

\vspace{1mm}
Now we assume that  $K'$ is not $W$-sound. Then there  is a cycle $K^*$ in $G^+$ containing  all of $W$ but not containing $x_4$. Consider $G'=G^+-x_4$.

\vspace{1mm} 
(i) Suppose $G'$ is $2$-connected. By Theorem \ref{thm3fh}, $G'\/$ has a $[v_4; x_3, w_1, w_2]\/$-$EPS\/$-graph $S' = E' \cup P'\/$ with $K^* \subseteq E'\/$. 
Set $E = E'\/$ and $P = P' \cup \{v_4x_4\}\/$ to obtain an $EPS\/$-graph of $G^+\/$ with $K^* \subseteq E\/$ and $v_4x_4\/$ is a pendant edge in $S\/$. Hence 
a hamiltonian cycle in $S^2\/$ with the required properties can be constructed. For $i=1,2$, note that if $w_ix_i\notin E(K^*)$, then $d_P(x_i)=0$ since 
$d_G(x_i)=2$ and $w_i\in K^*$. Observe also that $v_4x_i\in E(K^*)$ and $x_3x_i\in E(K^*)$ do not constitute any obstacle in this case.

\vspace{1mm} 
(ii) Suppose $G'\/$ is not $2\/$-connected. Let $B_y\/$ be the endblock in (the non-trivial block chain) $G'\/$ containing $K^*\/$, and let $t_4\/$ be the 
cutvertex of $G'\/$ belonging to $B_y\/$. Set $\widehat{G} = (G' - B_y) \cup \{u_4x_4\}\/$. Note that $\widehat{G}\/$ is  a non-trivial block chain and 
$\widehat{G} = ( G^+ - B_y) -x_4v_4\/$.

\vspace{1mm} 
Set $W^* = \{y, w_1, w_2, x_3, t_4\}\/$ and observe that $x_3 \not \in \{w_1, w_2\}\/$; otherwise, $G^+\/$ has a cycle containing $y, x_3, x_4\/$ 
(contradicting \textbf{E(1)}). In any case, $\widehat{G}\/$ has an $EPS\/$-graph $\widehat{S} = \widehat{E} \cup \widehat{P}\/$ with 
$d_{\widehat{P}}(t_4) \leq 1\/$ and $d_{\widehat{P}}(x_4) = 1\/$ by Lemma \ref{flelemma}(i).

\vspace{1mm} 
Now if $t_4\in\{w_1,w_2,x_3\}$, let $S_y=E_y\cup P_y$ be a $[t_4;r_4,s_4,y]$-$EPS$-graph of $B_y$ with $K^*\subseteq E_y$ where 
$\{r_4,s_4,t_4\}=\{w_1,w_2,x_3\}$, by Theorem \ref{thm3fh}.

\vspace{1mm} 
If, however,  $t_4\notin \{w_1,w_2,x_3\}$, we may assume without loss of generality that $K^*$ is $W^*$-sound (since $|W^*|=5$ and $K^*\supset W^*-t_4$). 
Consequently, let in this case $S_y = E_y \cup P_y\/$ be a $W^*$-$EPS\/$-graph of $B_y\/$ with $K^* \subseteq E_y\/$.

\vspace{1mm} 
In all cases, let an $EPS$-graph $S=E\cup P$ of $G^+$ be defined by $E=E_y\cup\widehat{E}$, $P = P_y \cup \widehat{P} \/$. We have $K^* \subseteq E\/$ and 
note that $d_P(w) \leq 1\/$ for every $w \in W^* - t_4\/$, and $d_P(t_4) \leq 2\/$ but  $d_P(t_4) \leq 1\/$ if $t_4 \in \{w_1, w_2, x_3\}\/$. It is now 
straightforward to see that in each of the cases in question, $S^2\/$ contains a hamiltonian cycle as required (see the argument at the end of case (i); 
moreover, $t_4x_i\in E(K^*)$ does not constitute an obstacle, $i=1,2$). This finishes case \textbf{(a1)}.

\vspace{2mm} 
Since the case $u_4\notin N(x_1)$ can be treated analogously, we are led to the following case.

\vspace{2mm} 
\hspace{5mm} {\bf (a2)} $u_4=w_1$ and $v_4=w_2$. Then $|V(K')|=6$. In view of case \textbf{(a1)}, we may assume that any cycle in $G^+$ containing 
$y, x_1, x_2, u_4, x_4, v_4\/$ has length  $6\/$.

\vspace{1mm} 
Suppose $H=G^+-x_1u_4$ is $2$-connected. Then $H$ has a cycle $C$ containing the edges $u_4x_4,x_4v_4,yx_1,yx_2,x_1w'_1\/$ (where $w'_1\neq u_4$). 
But this means that $|V(C)| > 6\/$ (because at least $2\/$ more edges are required to form the cycle $C\/$), a contradiction.

\vspace{1mm}
Thus $H$ is not 2-connected, and let $B_y\/$ and $B_4\/$ denote the endblocks of $H\/$ containing $y\/$ and $x_4\/$, respectively.

\vspace{1mm} 
Suppose $x_2$ is not a cutvertex of $H$. Since $\kappa(B_y)\geq 2$, it follows that $\{x_2,u_2,v_2\}\subset V(B_y)$. Now, we have a path $P=P(v_2,u_2)$ 
in $B_y\/$ with $x_2 \notin V(P)$. Since $d_G(x_4) =2\/$, $x_4 \not \in V(P)\/$; otherwise $u_4 \in V(P)\/$ as well and hence $x_4u_4 \in E(B_y \cap B_4)$ 
which is  impossible. Thus we obtain for $\{r_2, w_2\} = \{u_2,v_2\}\/$ a cycle $$K^* = (K' - w_2x_2) \cup P \cup \{r_2x_2\}$$  in $G^+\/$ containing $V(K')$ 
and $|V(K^*)| > 6\/$, contradicting the assumption at the beginning of this case. Thus $x_2\/$ is a cutvertex of $H\/$.

\vspace{1mm}
Observing that $d_{G^+}(x_2)=3\/$ and $\kappa(B_y) \geq 2\/$, we  conclude $d_{B_y}(x_2) =2\/$ and thus $x_2w_2 \in E(H) - E(B_y)\/$ is the other block of $H$ 
containing the cutvertex $x_2$. It now follows that $B_y\cap B_4=\emptyset$ since $x_2w_2 \not \in E(B_4)\/$. Without loss of generality $w_2 = v_2\/$; 
hence $u_2x_2 \in E(B_y)\/$.


\vspace{1mm} 
It now follows that $H -B_y\/$ is either a path of length $3\/$, or it is a block chain with $B_4\/$ being $2\/$-connected and $x_2v_2\/$ being a block.

\vspace{1mm} 
\hspace{8mm} {\bf (a2.1)} Suppose $x_3\in V(B_y)$. Let $K_y$ be a cycle in $B_y$ containing $y,x_1,x_2,x_3$ where we may assume that
$$K_y=yx_1w'_1\cdots x_3\cdots w'_2x_2y.$$ Note that $x_3=w'_1=w'_2$ is impossible because of $d_G(x_3)>2$.
If $x_3\neq w'_1$ and $x_3\neq w'_2$, then $B_y$ has an $[x_3;y,w'_1,w'_2]$-$EPS$-graph $S_y=E_y\cup P_y$ with $K_y\subseteq E_y$ by Theorem \ref{thm3fh}.
If $x_3=w'_1$ or $x_3=w'_2$, then $B_y$ has an $[x_3;y,w'_2]$-$EPS$-graph or an $[x_3;y,w'_1]$-$EPS$-graph $S_y=E_y\cup P_y$ with $K_y\subseteq E_y$ by 
Theorem \ref{thm2fh}, respectively. Likewise, if $d_G(u_4) >2\/$, then  $B_4\/$ has  a $[u_4; v_4]\/$-$EPS\/$-graph $S_4 = E_4 \cup P_4\/$ with 
$K^{(4)} \subseteq E_4\/$ where $K^{(4)}\/$  is a cycle in $B_4\/$ containing $u_4, x_4, v_4$, by Theorem \ref{thm1f}. If, however, $B_4$ is a bridge of $H$, 
then the path $P_4 = u_4x_4v_4\/$ has the only $EPS\/$-graph $S_4 = E_4 \cup P_4\/$ with $E_4 = \emptyset\/$.

\vspace{1mm} 
Setting $E=E_y\cup E_4$ and $P=P_y\cup P_4\cup\{x_1u_4\}$, we have an $EPS$-graph $S=E\cup P$ of $G^+$ with $d_P(x_1)=1$, $d_P(x_2)=d_P(x_3)=d_P(y)=0$, 
$d_P(w'_1) \leq 1\/$, $d_P(w'_2) \leq 1\/$,  $d_P(x_4) \in \{0,2\}\/$,  $d_P(u_4) \leq 2\/$ and $d_P(v_4) \leq 1\/$.  However, $d_P(x_4) = 2\/$ implies 
$d_P(v_4)=1\/$ and thus $x_4v_4\/$ is a pendant edge. Hence a hamiltonian cycle in $S^2\/$ with the required properties can be constructed.

\vspace{1mm} 
\hspace{8mm} {\bf (a2.2)} Suppose $x_3\in  V(B_4)$; thus $B_4$ is $2$-connected. Let $K_y$ be a cycle in $B_y$ containing $y,x_1,x_2$ where we may assume 
that $$K_y=yx_1w'_1\cdots w'_2x_2y.$$ Note that if $w'_1=w'_2$, then $d_{B_y}(w'_1)=2$. If $w'_1\neq w'_2$, then $B_y$ has an $[x_1;y,w'_1,w'_2]$-$EPS$-graph 
$S_y=E_y\cup P_y$ with $K_y\subseteq E_y$ by Theorem \ref{thm3fh}. If $w'_1=w'_2$, then $B_y$ has an $[x_1;y,w'_1]$-$EPS$-graph $S_y=E_y\cup P_y$ with 
$K_y\subseteq E_y$ by Theorem \ref{thm2fh}. Likewise, $B_4\/$ has a $[u_4; x_3, v_4]\/$-$EPS\/$-graph $S_4 = E_4 \cup P_4\/$ with $K^{(4)}\subseteq E_4$ 
where $K^{(4)}$ is a cycle in $B_4$ containing $x_3,u_4,x_4,v_4$. Setting $E=E_y\cup E_4$ and $P=P_y\cup P_4\cup \{x_1u_4\}$, we have an $EPS\/$-graph 
$S = E \cup P\/$ of $G^+\/$ and $S^2$ contains a hamiltonian cycle as required.

\vspace{2mm}
{\bf (b)}   $x_4 \in \{w_1, w_2\}\/$ but $w_1 \neq w_2\/$.

\vspace{1mm} 
Without loss of generality assume $x_4=w_1$ and hence $x_1x_4\in E(G)$ (the case $x_4=w_2$, $w_1\neq w_2$, can be solved by a symmetrical argument).
Note that $x_3=u_1=u_2$ cannot hold (see the argument in case \textbf{(1.1)(b)}).

\vspace{1mm} 
\hspace{5mm} {\bf (b1)} Suppose $v_4 \not \in N(x_2)\/$; i.e., $v_4 \neq w_2\/$. Let $K'\/$ be a cycle in $G^+\/$ containing $y, x_1, x_4, v_4, w_2, x_2\/$ in 
this order and let $W =\{y,x_4,v_4,w_2,x_3\}$. Then $K'$ is $W$-sound because of the supposition at the beginning of \textbf{(E)(1)}. By Theorem \ref{thm4fh}, 
$G^+$ has a $W$-$EPS$-graph $S=E\cup P$ with $K'\subseteq E\/$ and hence a hamiltonian cycle in $S^2\/$ with the required properties can be constructed.

\vspace{1mm} 
\hspace{5mm} {\bf (b2)} Suppose $v_4=w_2$. Assume first that $d_G(v_4)=2$. Let $K'$ be the cycle $yx_1x_4w_2x_2y\/$ and let $W = \{y, x_1, x_2, x_3, x_4\}\/$. 
Then $K'\/$ is $W\/$-sound. By Theorem \ref{thm4fh}, $G^+\/$ has a $W\/$-$EPS\/$-graph with $K' \subseteq E\/$.

\vspace{1mm}
Now assume that $d_G(v_4) > 2\/$. Let $z \in N(v_4) - \{x_4, x_2\}\/$.  There is a path $P(v_4, x_1)\/$ in $G\/$ from $v_4\/$ to $x_1\/$ via the vertex $z\/$ 
since $G\/$ is $2\/$-connected; $x_2 \not \in P(v_4, x_1)\/$ since $d_G(x_2) =2\/$.  Now $K^* = P(v_4, x_1) x_1yx_2v_4\/$ is a cycle in $G^+\/$ containing 
$N(x_4)\/$ but not $x_4\/$ itself. Hence $G'' = G^+ - x_4\/$ is $2\/$-connected.

\vspace{1mm} 
We may assume that $K^+$ is also a cycle in $G''$ containing $y,x_1,u_1,x_3,u_2,x_2$ in this order. If $x_3\neq u_1$ and $x_3\neq u_2$, then by Theorem 
\ref{thm2fh}, $G''$ has an $[x_3;u_1,u_2]$-$EPS$-graph $S''=E''\cup P''$ with $K^+\subseteq E''$. If $x_3=u_1$ or $x_3=u_2$, then by Theorem \ref{thm1f}, 
$G''$ has an $[x_3;u_2]$-$EPS$-graph or an $[x_3;u_1]$-$EPS$-graph $S''=E''\cup P''$ with $K^+\subseteq E''$, respectively.

\vspace{1mm}
Set $E=E''$ and $P = P'' \cup \{x_1x_4\}\/$. Then $S= E \cup P\/$ is an $EPS\/$-graph of $G^+\/$ such that $d_P(x_1) =1\/$,  $ d_P(x_3)=0= d_P(y)\/$ and 
$d_P(w) \leq 1\/$ for $w \in \{x_2, u_1, u_2\}\/$ and $x_1x_4\/$ is a pendant edge in $S\/$. In either case,  a hamiltonian cycle in $S^2\/$ with the required 
properties can be constructed.

\vspace{2mm}
{\bf (c)}   $N(x_4) = \{x_1, x_2\}\/$.

\vspace{1mm} 
Clearly $G'' = G^+ - x_4\/$ is $2\/$-connected. Let $K''\/$ be a cycle in $G''\/$ containing $y, x_1, x_2, x_3\/$, and let $u_1 \in V(K'') \cap N_G(x_1)\/$. 
Without loss of generality, $u_1 \neq x_3\/$: for $d_G(x_3) > 2\/$ implies $\{x_1, x_2\} \not \subset N(x_3)\/$.

\vspace{1mm} 
Then $G''\/$ has an $[x_3;u_1]$-$EPS$-graph $S''=E''\cup P''$ with $K''\subseteq E''$. Set $E = E''\/$ and $P = P'' \cup \{x_1x_4\}\/$. Then $S=E \cup P\/$ is 
an $EPS\/$-graph of $G^+\/$ with $d_P(y) =0 = d_P(x_3) = d_P(x_2)\/$ and $x_1x_4\/$ being a pendant edge in $S\/$. Hence a hamiltonian cycle in $S^2\/$ with 
the required properties can be constructed.

\vspace{3.88mm}
{\bf (1.3)}  $d_G(x_3)=2, \ d_G(x_4) =2\/$.

\vspace{2mm}
Recall that $x_3,x_4$ are not on the same cycle containing $y, x_1, x_2\/$. For each $i=3, 4 \/$, let $l_i\/$ denote the length of a longest cycle in $G^+\/$ 
containing $y, x_1, x_2, x_i\/$.

\vspace{2mm}
{\bf (a)} Suppose $l_3\geq 7$ or $l_4\geq 7$; without loss generality assume that $l_3\geq 7$. Recall that $$K^+=yx_1u_1 \cdots u_3x_3v_3\cdots u_2x_2y$$
Then either $u_1 \not \in \{u_3, x_3\}\/$ or $u_2 \not \in \{v_3, x_3\}\/$. Without loss of generality, assume that $u_1 \not \in \{u_3, x_3\}\/$.

\vspace{2mm} 
\hspace{5mm} {\bf (a1)} Assume that $G' = G^+ -x_4\/$ is $2\/$-connected.

\vspace{2mm} 
Set $W=\{y,u_1,u_2,u_3,q_4\}$, where $q_4\in \{u_4,v_4\}$. Note that $|\{y,u_1,u_2,u_3\}|=4$. 

\vspace{2mm} 
Suppose $q_4$ exists such that $|W|=4$, say for $q_4=u_4$. Then $u_4\in\{u_1,u_2,u_3\}$ and $G'$ has a $[u_4;w_1,w_2]$-$EPS$-graph $S'=E'\cup P'$ with 
$K^+\subseteq E'$, where $\{u_4,w_1,w_2\}=\{u_1,u_2,u_3\}$, by Theorem \ref{thm2fh}.

\vspace{2mm} 
Now suppose that $|W|=5$ and $K^+$ is $W$-sound in $G'$ for some choice of $q_4$, say for $q_4=u_4$. Then by Theorem \ref{thm4fh} there is a $W$-$EPS$-graph 
$S'=E'\cup P'$ of $G'$ with $K^+\subseteq E'$.

\vspace{2mm} 
In both cases, taking $E=E'$ and $P=P'\cup\{x_4u_4\}$, we have an $EPS$-graph $S=E\cup P$ of $G^+$ such that $d_P(w)\leq 1$ for all $w\in W-\{u_4\}$,  
$d_P(x_4)=1$ and $d_P(u_4)\leq 2$. Hence a required hamiltonian cycle in $S^2$ can be constructed; it can be made to contain $x_4u_4$ and $u_3x_3$.


\vspace{2mm}
Hence we assume that $|W|=5$ and $K^+$ is not $W$-sound in $G'$ for any choice of $q_4\in\{u_4,v_4\}$. Then there is another cycle $K'$ in $G'$ such that 
$V(K')\supseteq W$. We may assume that $q_4=u_4$ and $x_3\notin K'$. Then by Theorem \ref{thm3fh}, $G'$ contains a $[u_3;u_1,u_2,u_4]$-$EPS$-graph 
$S'=E'\cup P'$ with $K'\subseteq E'$. Taking $E=E'$ and $P=P'\cup\{x_4u_4\}$, we have an $EPS$-graph $S=E\cup P$ of $G^+$. Note that $x_4$ is a pendant vertex 
in $S$ and either $x_3$ is a vertex in $E$, or else it is a pendant vertex in $S$. Hence a required hamiltonian cycle in $S^2$ can be constructed. For 
$i=1,2$, also note that if $u_ix_i\notin E(K')$, then $d_P(x_i)=0$ since $d_G(x_i)=2$ and $u_i\in K'$. Observe also that $u_3x_i\in E(K')$ and 
$u_4x_i\in E(K')$ do not constitute any obstacle in this case.



\vspace{2mm} 
\hspace{5mm} {\bf (a2)} Assume that $G' = G^+ -x_4\/$ is not $2\/$-connected.

\vspace{2mm} In view of case \textbf{(a1)}, we may assume, by symmetry, that $G^+ -x_3\/$ is also not $2\/$-connected.

\vspace{2mm}
Let $K^{(i)}$ denote a cycle containing $y, x_1, x_2, x_i\/$ where $i \in \{3,4\}\/$. Let $B_i\/$ be the block of $G^+ - x_i\/$ with $K^{(7-i)} \subset B_i$.  
Let  $G_i, \ G_i'\/$ denote the block chains in $G^+ - x_i -B_i\/$ (possibly $G_i = \emptyset\/$ or $G_i' = \emptyset\/$) which contain $\{u_i, c_i\}\/$ and 
$\{v_i, c_i'\}$ respectively,  where $c_i,c_i'$ denote the cutvertices of $G^+ - x_i$ belonging to $B_i$, provided $G_i\neq \emptyset$, 
$G_i' \neq \emptyset\/$. If $ G_i = \emptyset\/$, then $u_i = c_i\/$ and is not a cutvertex, and likewise $v_i = c_i'\/$ if $G_i' = \emptyset\/$.

\vspace{2mm}
We observe that $K^{(7-i)}\/$ is edge-disjoint from $G_i \cup G_i'\/$, $i=3,4\/$  and that  $G_3 \cup G_3'\/$ and $G_4 \cup G_4'\/$ are edge-disjoint (since 
every block of $G_i \cup G_i'$ contains an edge of $K^{(i)}$). Finally, if $C_y$ (in $G^+$) is a cycle containing $y$, then 
$E(C_y\cap (G_i\cup G_i')) = \emptyset \/$ for at least one $i \in \{3, 4\}\/$; otherwise, $C_y \supset \{x_3, x_4\}\/$, contrary to \textbf{(E)(1)}. 
Without loss of generality $C_y$ is one of the cycles $K^{(3)}$, and we may also assume that $K^{(3)}=K^+$ (see the beginning of \textbf{(a)}).

\vspace{2mm} 
Set $W=\{y,u_1,u_2,u_3,x_4\}$. The definition of $W$ together with the last sentences of the preceding paragraph ensure that $|W|=5$ and $K^{(3)}=K^+$ is 
$W$-sound in $G^+$.

\vspace{2mm}
Set $G_0 = G_4 \cup G_4' \cup \{u_4x_4v_4\}\/$; it is a block chain.

\vspace{2mm}  
\hspace{8mm} {\bf (a2.1)} Suppose $G_0 \/$  is a path with   $3 \leq l(G_0) \leq 4\/$.

\vspace{1mm} 
Then by Theorem \ref{thm4fh}, $G^+\/$ has a $W\/$-$EPS\/$-graph $S = E \cup P\/$ with $K^{(3)} \subseteq E\/$ and $d_P(x_4) \leq 1\/$. If $d_P(x_4)=0\/$, 
then $x_4\/$ is in $E\/$, and one of its neighbors is $2\/$-valent because  $l(G_0) \geq 3\/$. If $d_p(x_4)=1\/$, then $x_4\/$ is a pendant vertex in $S\/$. 
In either case, a required hamiltonian cycle in $(G^+)^2\/$ can be constructed.

\vspace{2mm} 
\hspace{8mm}{\bf (a2.2)} Suppose $G_0\/$ is a path with  $l(G_0) \geq 5\/$, or $G_0\/$ is a block chain having non-trivial blocks.

\vspace{1mm} 
Replace $G_0\/$ in $G^+\/$ by a path $P_4 = c_4u_4x_4v_4c_4'\/$ to obtain the graph $G^*$ (note that $|E(G_0)|\geq 5$). Observe that the cycle $K^{(3)}\/$ 
in  $G^*\/$ passes through  the vertices $y, x_1, x_2, x_3\/$. Then as in case \textbf{(a2.1)}, $G^*\/$ has a $W\/$-$EPS\/$-graph $S^* = E^* \cup P^*\/$ 
with $d_{P^*}(x_4)=0\/$ or $d_{P^*}(x_4) =1\/$.

 \vspace{1mm}
 (i) If $d_{P^*}(x_4)=0\/$, then  $P_4 \subset E^*\/$. Since $G_0\/$ is a block chain, by Lemma \ref{flelemma}(ii), $G_0\/$ contains a $JEPS\/$-graph $S_0 = J_0 \cup E_0 \cup P_0\/$ such that $d_{P_0}(c_4)=0= d_{P_0}(c_4')\/$. Moreover, in constructing $S_0\/$ by proceeding block by block, one can achieve $d_{P_0}(u_4) \leq 1\/$, $d_{P_0}(v_4) \leq 1\/$.     In this case, we obtain a $W\/$-$EPS\/$-graph $S = E \cup P\/$ of $G^+\/$ by setting $E = (E^* -P_4) \cup J_0 \cup E_0 \/$ and $P = P^* \cup P_0\/$. Here $d_{P}(x_4) =0\/$, $d_P(u_4) \leq 1, d_P(v_4) \leq 1, d_P(c_4) \leq 2, d_P(c_4') \leq 2\/$ and a required hamiltonian cycle in $(G^+)^2\/$ can be constructed.

\vspace{2mm}
(ii) If $d_{P^*}(x_4)=1\/$, then  $V(P_4) \subseteq V(P^*) \/$. Hence either $u_4x_4 \not \in E(P^*)\/$ or $v_4x_4 \not \in E(P^*)\/$. Suppose $v_4x_4 \not \in E(P^*)\/$ (so that $u_4x_4 \in E(P^*)\/$).  By Lemma \ref{flelemma}(i), $G_4 \cup \{u_4x_4\}\/$ (respectively $G_4' \/$) has an $EPS\/$-graph $S^{(4)} = E^{(4)} \cup P^{(4)}\/$ (respectively $S'^{(4)} = E'^{(4)} \cup P'^{(4)}\/$) such that $d_{P^{(4)}}(c_4) \leq 1, \ d_{P^{(4)}}(u_4) \leq 2, \ d_{P^{(4)}}(x_4) = 1\/$ with $u_4x_4\/$ being a pendant edge in $S^{(4)}\/$ and $d_{P'^{(4)}}(c_4') \leq 1, \ d_{P'^{(4)}}(v_4) \leq 1\/$.  Now, if we take $E =  E^* \cup E^{(4)} \cup E'^{(4)} \/$ and $P = (P^* -\{u_4, v_4\}) \cup P^{(4)} \cup P'^{(4)}\/$, we have an $EPS\/$-graph $S= E \cup P\/$ of $G^+\/$ with $d_{P} (w) \leq 1\/$ for every $w \in W\/$ from which a required hamiltonian cycle in $(G^+)^2\/$ can be constructed (take note that $c_4u_4, v_4c_4' \in E(P^*)\/$ resulting in $d_P(c_4) \leq 2\/$ and $d_P(c_4') \leq 2\/$; and $d_P(x_i) \leq 1\/$ is guaranteed by the assumption $d_G(x_i) =2\/$, $i=1, 2\/$).

\vspace{2mm} 
In view of case \textbf{(1.3)(a)} solved, we may assume from now on that $l_3\leq l_4$ and hence we are left with the following case.

\vspace{3.88mm}
{\bf (b)}  Suppose $4 \leq l_3 \leq l_4 \leq 6\/$.

\vspace{2mm} 
\hspace{5mm}{\bf (b1)}  Suppose $l_3=6$.

\vspace{3mm} \hspace{8mm}
{\bf (b1.1)} Suppose $u_3=u_4=u_1$ and $v_3=v_4=v_2$. Set  $G^* = G -x_4\/$.

\vspace{2mm}
$\kappa(G^*) =2\/$ since  $N(x_4) = N(x_3)\/$. By induction, $G^*\/$ has the ${\cal F}_4$-property; that is, there exists an $x_1x_2\/$-hamiltonian path 
$P(x_1, x_2)\/$ in $(G^*)^2\/$ containing different edges $x_3z_3, u_4z_4 \in E(G^*)\/$. We may write $$P(x_1, x_2) = x_1 \cdots st  \cdots x_2$$  where 
$\{s, t\} = \{ u_4, z_4\}.$ Then $$x_1 \cdots sx_4t \cdots x_2$$ is a required hamiltonian path in $G^2$; it contains $x_3z_3$ because $P(x_1, x_2)$ does.

\vspace{3mm} 
\hspace{5mm}{\bf (b1.2)} Suppose $u_3=u_1, \ v_3= v_2\/$ and  $u_4=v_1, \ v_4 = u_2\/$.

\vspace{2mm}
Consider $G^- = G^+ - \{x_1v_1, x_2u_2\}\/$. If there is a path $P(s,t)\/$ from $s \in \{v_1, u_2\}\/$ to $t \in \{u_1, v_2\}\/$ in $G^-\/$, then either 
$l_3 > 6\/$ or $l_4 >6\/$, or $G^+\/$ has a cycle containing both $x_3\/$ and $x_4\/$. Thus $x_3\/$ and $x_4\/$ belong to different components of $G^-\/$.
Let $G_i\/$ denote the  component of $G^-\/$ containing the vertices $u_i,  x_i, v_i\/$, $i \in \{3,4\}\/$. We reach the same conclusion when considering 
$G^+ - \{x_1u_1, x_2v_2\}\/$ instead of $G^-\/$. Since $N_G(x_1) \not \subseteq V_2(G)\/$, $N_G(x_2) \not \subseteq V_2(G)\/$, we may assume without loss of 
generality that $d(v_1) >2\/$ or $d(u_2) > 2$ (otherwise, $x_3$ and $x_4$ switch their roles) and hence both $v_1, u_2\/$ are not $2\/$-valent (otherwise, 
$v_1\/$ or $u_2\/$ would be a cutvertex of $G\/$). It follows that $G_4\/$ is $2\/$-connected. Likewise, $G_3$ is also 2-connected.

\vspace{1mm} 
There is a cycle $C^{(4)}$ in $G_4$ containing $u_4,x_4,v_4$ and there is a cycle $C^{(3)}$ in $G_3$ containing $y,x_1,x_2,u_3,x_3,v_3$. 
By Theorem \ref{thm1f}, $G_i$ has a $[u_i;v_i]$-$EPS$-graph $S_i=E_i\cup P_i$  with $C^{(i)}\subseteq E_i$, $i=3,4$. Note that $d_{P_3}(z)=0$ 
for $z\in \{y,x_1,x_2, x_3\}$.

\vspace{1mm}  
Now set $E=E_3\cup E_4$ and $P=P_3\cup P_4\cup \{x_1v_1\}$. Then $S=E \cup P\/$ is an $EPS\/$-graph of $G^+\/$ with $C^{(3)} \cup C^{(4)} \subseteq E\/$ 
and a required hamiltonian cycle in $(G^+)^2\/$ containing $x_4v_1, x_3v_2\/$ can be constructed.

\vspace{3mm} 
\hspace{5mm}{\bf (b1.3)} Suppose $u_3=u_1=u_4,v_3=v_2$ and $v_4=u_2$ (the case $u_3=u_1$, $u_4=v_1$ and $v_2=v_3=v_4$ is symmetric).

\vspace{2mm} 
This subcase is impossible;  otherwise, it gives rise to a cycle containing $y, x_3, x_4\/$, a contradiction to  the assumption (just consider in $G\/$ a path 
from $x_1\/$ to $u_2\/$ avoiding $u_1\/$).

\vspace{1mm} 
It is straightforward to see that $x_i\notin N(x_j)$ for  $i=3,4$ and $j=1,2$ for all choices of $i$ and $j$; otherwise, $l_i>6$ or there exists a cycle
containing $y,x_3,x_4$. Therefore, subcase \textbf{(b1)} is finished.

\vspace{3mm}
{\bf (b2)}  Suppose $l_3=5$.

\vspace{2mm}   
We may assume without loss of generality that $u_3 = x_1, \ x_3 = u_1 \/$ and $v_3 = v_2\/$. 

\vspace{2mm}
Suppose $d_G(v_3)=2$. Consider $G'=G-\{x_3,v_3\}$; it is a non-trivial block chain with pendant edges $x_1v_1, x_2u_2$. By Corollary \ref{flecor}(ii), there 
exists a hamiltonian path  $P(x_1, x_2) \subseteq (G')^2\/$ starting with $x_1v_1\/$ and ending with $u_2x_2\/$. We proceed block by block to construct 
$P(x_1, x_2)\/$ such that $x_4z_4 \in E(G) \cap P(x_1, x_2)\/$ and $x_4z_4 \not \in \{x_1v_1, u_2x_2\}\/$: this is clear if $x_4\/$ is a cutvertex of $G'\/$; 
and if $x_4 \in V(B_4)\/$ where $B_4 \subseteq G'\/$ is a $2\/$-connected block containing the cutvertices $c_4, c_4'\/$ of $G'\/$, one uses a hamiltonian 
path $P(c_4, c_4')\/$ in  $(B_4)^2\/$ containing an edge incident to $x_4\/$ (Theorem \ref{fletheorem4}(i)). Then $$(P(x_1, x_2)-u_2x_2) u_2v_3x_3x_2$$ is a 
required hamiltonian path in $G^2$.

\vspace{1mm}  
If $d_G(v_3)>2$, then $G^{(0)}=G-\{x_1,x_2,x_3\}$ is connected (or else $v_3$ is a cutvertex of $G\/$). Any $v_3u_2\/$-path $P(v_3, u_2) \subset G^{(0)}\/$ 
can be extended to a cycle $yx_1x_3P(v_3, u_2)x_2y\/$ of length $\geq 6$, contradicting the assumption of this subcase.

\vspace{3mm}
{\bf (b3)}  Suppose $l_3 =4\/$.

\vspace{2mm}
In this case, let $G'=G-x_3$. Operating with $P(x_1,x_2)\subseteq (G')^2$ as in case \textbf{(b2)}, we obtain an ${\cal F}_4\/$ $x_1x_2\/$-hamiltonian path 
$(P(x_1,x_2) - u_2x_2) u_2x_3x_2\/$ in $G^2$.

\vspace{3mm}  
{\bf (1.4)} $d_G(x_3)=2, d_G(x_4)>2$.

\vspace{1mm} 
This case is symmetrical to the case \textbf{(1.2)}.

\vspace{5mm} 
{\bf (E)(2)} Suppose $x_3$ and $x_4$ are in $K^+$.

\vspace{2mm} 
Without loss of generality, assume that $$K^+ = yx_1u_1 \cdots z_3x_3 \cdots x_4z_4 \cdots u_2x_2y.$$ As for the definition of $x_3^*, x_4^*$ see the 
paragraph preceding the statement of Lemma \ref{lemma4cd}.

\vspace{2mm}
{\bf (2.1)} $x_3 \neq u_1$ and $x_4 \neq u_2$.

\vspace{2mm}  
{\bf (a)} Suppose either $u_{i-2}\not\in N_G(x_i)$, or $u_{i-2}\in\in N_G(x_i)$ and $d_G(x_i)>2$ for some $i \in \{3,4\}$. Without loss of generality, 
assume that $i=4$.

\vspace{1mm}
If $u_1\neq x_3^*$, set $W=\{y,u_1,u_2,x_3^*,x_4^*\}$. Then $|W|=5$ and $K^+$ is $W$-sound, so by Theorem \ref{thm4fh}, $G^+$ has a $W$-$EPS$-graph 
$S=E\cup P$ with $K^+ \subseteq E$.

\vspace{1mm}
If $u_1=x_3^*$, then $d_G(x_3)=2$ since $x_3\neq u_1$ by supposition. Now, let $S=E\cup P$ be an $[x_4^*;u_1,u_2]$-$EPS$-graph of $G^+$ with $K^+\subseteq E$
by Theorem \ref{thm2fh}.

\vspace{1mm} 
In either case, a required hamiltonian cycle in $(G^+)^2$ can be constructed.

\vspace{2mm} 
{\bf (b)} Suppose $u_{i-2}\in N_G(x_i)$ and $d_G(x_i)=2$ for $i=3,4$.

\vspace{2mm} 
If $w_4$ is the predecessor of $x_4$ in $K^+$ and $w_4\neq x_3$, then let $S= E\cup P$ be a $[x_1 ; u_1, u_2, w_4]$-$EPS$-graph with $K^+\subseteq E$ 
by Theorem \ref{thm3fh}. If $w_4=x_3$, then let $S=E \cup P$ be an $[x_1;u_1,u_2]$-$EPS$-graph with $K^+\subseteq E$ by Theorem \ref{thm2fh}. Hence 
a required hamiltonian cycle in $(G^+)^2$ can be constructed from $S$.

\vspace{2mm}
{\bf (2.2)} $x_3=u_1$ and $x_4\neq u_2$.

\vspace{2mm}
{\bf (a)} Suppose either $u_2\notin N_G(x_4)$, or $u_2\in N_G(x_4)$ and $d_G(x_4)>2$.

\vspace{1mm} 
\hspace{5mm} {\bf (a1)} $x_3x_4 \in E(G)$.

\vspace{1mm} 
If $d_G(x_4)>2$, then $d_G(x_3)=2$ and we choose an $[x_1;x_4,u_2]$-$EPS$-graph $S=E\cup P$ of $G^+$ with $K^+\subseteq E$ by Theorem \ref{thm2fh}. If, 
however $d_G(x_4)=2$, we choose an $[x_1;x_3,z_4,u_2]$-$EPS$-graph $S=E\cup P$ of $G^+$ with $K^+\subseteq E$ by Theorem \ref{thm3fh}. In either case, 
$S^2$ contains a required hamiltonian cycle.

\vspace{2mm} 
\hspace{5mm} {\bf (a2)} $x_3x_4\notin E(G)$.

\vspace{1mm} 
Here $w_3$ is the successor of $x_3$ in $K^+$. Let $S=E\cup P$ be an $[x_3;u_2,w_3,x_4^*]$-$EPS$-graph with $K^+\subseteq E$ by Theorem \ref{thm3fh}. Also 
here, $S^2$ contains a required hamiltonian cycle; it contains $x_3v\in E(G)$ which is consecutive to $x_1x_3$ in the eulerian trail of the component of $E$ 
containing $K^+$ (possibly $v=w_3$) and it contains $x_4z_4$.

\vspace{2mm} 
{\bf (b)} Suppose $u_2\in N(x_4)$ and $d_G(x_4)=2$.

\vspace{1mm} 
\hspace{5mm} {\bf (b1)} $x_3x_4\in E(G)$.

\vspace{1mm} 
Let $H = G-x_4\/$. Suppose $H\/$ is $2\/$-connected. Then by induction, $H\/$ has an ${\cal F} _4\/$ $x_1x_2\/$-hamiltonian path $P(x_1, x_2)\/$ in $H^2\/$ 
containing $x_3w_3\/$ and $u_2w_2\/$ which are edges of $H\/$. By deleting $u_2w_2\/$ from $P(x_1, x_2)\/$ and joining $x_4\/$ to $u_2, w_2\/$, we obtain 
an ${\cal F} _4\/$ $x_1x_2\/$-hamiltonian path in $G^2\/$ containing $x_3w_3, x_4u_2\/$ which are edges of $G\/$.

\vspace{1mm} 
Suppose $H\/$ is not $2\/$-connected. Then $H\/$ is a non-trivial block chain with endblock $B_i\/$ containing $u_i\/$; $u_i\/$ is not a cutvetex of $H\/$, 
$i=1, 2\/$. Let $c_i\/$ denote the cutvertex of $H\/$ which is contained in $B_i\/$, $i = 1,  2\/$. Set $B_{1,2} = H - (B_1 \cup B_2)\/$. If $c_1 = c_2\/$, 
then set $B_{1,2}=c_1\/$. In any case, $c_1\/$ and  $c_2\/$ are not cutvertices of $B_{1,2}\/$.

\vspace{1mm} 
By supposing $x_i \neq c_i\/$ (and thus $B_i\/$ is $2\/$-connected) we apply Theorem  \ref{fletheorem4} to conclude that $(B_i)^2\/$ has an ${\cal F} _3\/$ 
$x_ic_i\/$-hamiltonian path $P(x_i,c_i)\/$,   $i=1,2\/$ containing $x_3w_3, u_2w_2\/$ respectively,  which are edges of $G\/$. Let $P(c_1, c_2)\/$ denote 
a $c_1c_2\/$-hamiltonian path in $(B_{1,2})^2\/$. By deleting the edge $u_2w_2\/$ from the $x_1x_2\/$-hamiltonian path $P(x_1,c_1)P(c_1,c_2)P(x_2,c_2)\/$  
in $(G-x_4)^2\/$ and joining $x_4\/$ to $u_2, w_2\/$, we obtain an ${\cal F} _4\/$ $x_1x_2\/$-hamiltonian path in $G^2\/$ containing $x_3w_3, x_4u_2\/$ 
which are edges of $G$. Now suppose $x_1=c_1$ or $x_2=c_2$; i.e., $d_G(u_1) = 2$ or $d_G(u_2) =2\/$. In this case we consider $G^+\/$ and choose 
an $[x_1; u_1, u_2]$-$EPS$-graph $S=E\cup P$ of $G^+$ with $K^+\subseteq E\/$ by Theorem \ref{thm2fh}. Hence $S^2$ contains a hamiltonian cycle as required.

\vspace{2mm}
\hspace{5mm} {\bf (b2)} $x_3x_4 \notin E(G)$.

\vspace{1mm} 
If $w_3\neq w_4$, then we set $S=E\cup P$ to be an $[x_3;u_2,w_3,w_4]$-$EPS$-graph of $G^+$ with $K^+\subseteq E$ by Theorem \ref{thm3fh}. If $w_3=w_4$, then 
we set $S=E\cup P$ to be an $[x_3;u_2,w_3]$-$EPS$-graph of $G^+$ with $K^+\subseteq E$ by Theorem \ref{thm2fh}. Here $w_3$ is the successor of $x_3$ and 
$w_4$ is the predecessor of $x_4$ in $K^+$. Hence $S^2$ yields a required hamiltonian cycle unless $w_3=w_4$ and $d_G(w_3)>2$, in which case $d_G(x_3)=2$ 
holds, and we operate with an $[x_1; w_3, u_2]$-$EPS$-graph by Theorem \ref{thm2fh}. This settles case \textbf{(2.2)}.

\vspace{2mm} 
Since the case $x_3\neq u_1$ and $x_4=u_2$ is symmetrical to the case \textbf{(2.2)} just dealt with, we are left with the following case.

\vspace{3mm}
{\bf (2.3)} $x_3=u_1$ and $x_4=u_2$.

\vspace{2mm} 
\hspace{5mm} {\bf (a)} $d_G(x_3)=2$.

\vspace{2mm} 
\hspace{8mm} {\bf (a1)} $x_3x_4 \notin E(G)$.

\vspace{2mm} 
Choose an $[x_4;u_3,u_4]$-$EPS$-graph $S=E\cup P$ of $G^+$ with $K^+\subseteq E$ by Theorem \ref{thm2fh} if $u_3\neq u_4$, and an $[x_4;u_3,u_4]$-$EPS$-graph 
$S=E\cup P$ of $G^+$ with $K^+\subseteq E$ by Theorem~\ref{thm1f} if $u_3=u_4$; here $u_3$ is taken to be the successor of $x_3$ and $u_4$ the predecessor 
of $x_4$ in $K^+$. Then $S^2$ yields a required hamiltonian cycle unless $u_3=u_4$ and $d_G(u_3)>2$. In this case $d_G(x_4)=2$ and we may operate with 
an $[x_2;u_3]$-$EPS$-graph to obtain a required hamiltonian cycle in $S^2$ by Theorem \ref{thm1f}.

\vspace{2mm} 
\hspace{8mm} {\bf (a2)} $x_3x_4\in E(G)$.

\vspace{1mm} (i) Suppose $d_G(x_4)>2$.

\vspace{1mm}   
$G-x_3\/$ is a block chain in which $x_1\/$ and $x_4\/$ are not cutvertices and belong to different endblocks. However, the endblock containing $x_4\/$ is 
$2\/$-connected since $d_G(x_4) >2\/$; and it contains $x_2\/$ as well which is not a cutvertex of $G-x_3\/$ either. Therefore, $G^+-x_3$ is $2$-connected. 
Set $$H = (G^+ - \{y, x_1, x_3\}) \cup \{x, xv_1, xx_2\}. $$  $H\/$ is $2\/$-connected since $G^+-x_3\/$ is $2\/$-connected. By Theorem \ref{fletheorem3}, 
$H^2\/$ has a hamiltonian cycle $C\/$ containing $v_1x, xx_2, x_4w_4\/$ which are edges of $H\/$. Now $(C- x) \cup \{v_1x_3x_1yx_2\} \/$ is a hamiltonian 
cycle in $(G^+)^2\/$ with the required properties.

\vspace{1mm}
(ii) Suppose $d_G(x_4)=2$.

\vspace{1mm}
Let $H\/$ be the graph obtained from $G^+\/$ by deleting $y, x_2, x_3, x_4\/$.  Then $H\/$ is a non-trivial block chain containing $x_1\/$ which is not 
a cutvertex of $H\/$.  By Corollary \ref{flecor}(i), $H^2\/$ has a hamiltonian cycle $C\/$ containing the edge $x_1v_1\/$ (which is an edge of $G\/$). This 
implies that the cycle $yx_1(C - x_1v_1) v_1x_3x_4x_2y\/$ is a hamiltonian cycle in $(G^+)^2\/$ having the required properties.

\vspace{3mm}
{\bf (b)} $d_G(x_3)>2$, hence $d_G(x_4)>2$; otherwise we are back to {\bf (a)} above, by symmetry. Then $x_3x_4\notin E(G)$.

\vspace{2mm}
Suppose $G'=G-x_1$ is $2$-connected. Then by induction, $G'$ has an ${\cal F}_4$ $v_1x_2$-hamiltonian path $P(v_1, x_2)$ in $(G')^2$ containing $x_3w_3$ 
and $x_4w_4$ which are edges of $G'$. Now $\{x_1v_1\}\cup P(v_1, x_2)$ is an ${\cal F} _4\/$ $x_1x_2$-hamiltonian path in $G^2$ containing $x_3w_3, x_4w_4$ 
which are edges of $G$.

\vspace{2mm}
Now suppose $G'=G-x_1$ is not $2$-connected. Then $G'$ is a non-trivial block chain with $x_3,v_1$ in different endblocks and not cutvertices. 
Note that the block containing $x_3$ is 2-connected and at least one block contaning $x_4$ is 2-connected, since $d_G(x_3)>2$ and $d_G(x_4)>2$.

\vspace{2mm}
\hspace{5mm} {\bf (b1)} Suppose $x_2$ is a cutvertex of $G'$. Let $G_1$ and $G_2$ be the components of $G'-x_2$ with either $x_3,x_4\in V(G_1)$ and
$v_2,v_1\in V(G_2)$, or $x_3,v_2\in V(G_1)$ and $x_4,v_1\in V(G_2)$ (note $d_{G'}(x_2)=2$). Observe that in the first case $v_2=v_1$ is possible. 
However, $v_1=x_4$ is impossible because of the assumptions of this case \textbf{(b)}; i.e., $d_G(x_4)>2$. By the same token $v_2=x=3$ is impossible.

\vspace{2mm}
Suppose $x_3,x_4\in V(G_1)$ and $v_2,v_1\in V(G_2)$. Then by Theorem \ref{fletheorem4}(ii) or Corollary~\ref{flecor}(ii), respectively, $(G_1)^2$ has an 
$x_3x_4$-hamiltonian path $P_1$ containing an edge $x_3w_3\in E(G)$. If $G_2=K_1=v_1$, then we set $P=P_1\cup\{x_2x_4,x_3v_1,v_1x_1\}$. If $G_2=K_2=v_2v_1$, 
then we set $P=P_1\cup\{x_2x_4,x_3v_1,v_1v_2,v_2x_1\}$. Otherwise, by Theorem~\ref{fletheorem3} or Corollary \ref{flecor}(i), respectively, $(G_2)^2$ has 
a hamiltonian cycle $C_2$ containing an edge $t_1v_1\in E(G)$. Then we set $P=P_1\cup C_2\cup\{x_2x_4,x_3v_1,t_1x_1\}-\{t_1v_1\}$. In all cases $P$ is an 
${\cal F}_4$ $x_1x_2$-hamiltonian path in $G^2$ containing $x_3w_3, x_4x_2$ which are edges of $G$ as required.

\vspace{2mm}
Suppose $x_3,v_2\in V(G_1)$ and $x_4,v_1\in V(G_2)$. Then we apply an analogous strategy as in the preceding case using Theorems~\ref{fletheorem3}, 
\ref{fletheorem4} and Corollary~\ref{flecor}, but considering $G_1$ instead of $G_2$ and vice versa.

\vspace{2mm}
\hspace{5mm} {\bf (b2)} Suppose $x_2$ is not a cutvertex of $G'$. Let $B_2$ be the 2-connected block containing $x_2$.

\vspace{2mm}
(i) Suppose $x_3\in V(B_2)$. Let $t$ be the cutvertex of $G'$ in $B_2$; possibly $t=x_4$, $t\notin \{x_2,x_3\}$ in any case. We define the block chain $G_1$ 
such that $G'=B_2\cup G_1$ and $B_2\cap G_1=\{t\}$. If $t=x_4$, then $(B_2)^2$ has an $x_2t$-hamiltonian path $P_2$ containing $x_3w_3\in E(G)$ by 
Theorem~\ref{fletheorem4}(i). If $t\neq x_4$, then by induction $(B_2)^2$ has an $x_2t$-hamiltonian path $P_2$ containing $x_3w_3, x_4w_4$ which are different 
edges of $G$. In both cases $(G_1)^2$ has a $tv_1$-hamiltonian path starting with $tw\in E(G)$, by Theorem~\ref{fletheorem4}(ii) or 
Corollary~\ref{flecor}(ii), respectively. Then $P=P_2\cup P_1\cup \{v_1x_1\}$ is an ${\cal F}_4$ $x_1x_2$-hamiltonian path in $G^2$ containing 
$x_3w_3, x_4w_4$ which are edges of $G$ as required. Note that if $t=x_4$, then $x_4w_4=tw$.

\vspace{2mm}
(ii) Suppose $x_3\notin V(B_2)$. If $B_2$ is not an endblock, then $t,t'$ denote the cutvertices of $G'$ in $B_2$ and we define block chains $G_0$, $G_1$ 
such that $G'=G_1\cup B_2\cup G_0$, $x_3\in V(G_1), v_1\in V(G_0)$ and $G_1\cap B_2=t$, $B_2\cap G_0=t'$. If $B_2$ is an endblock, then we proceed analogously: we set $G_0=\emptyset$ and $t'=v_1$ in this case. Note that $t=x_4$ ot $t'=x_4$ is possible.

\vspace{2mm}
If $t'\neq x_4$, then by Theorem \ref{fletheorem4}(i) $(B_2)^2$ has an $x_2t$-hamiltonian path $P_2$ containing $t'w'\in E(G)$ for $t=x_4$ and by induction 
$(B_2)^2$ has an ${\cal F}_4$ $x_2t$-hamiltonian path $P_2$ containing $t'w', x_4w_4$ which are different edges of $G$ for $t\neq x_4$. By the same token 
$(G_1)^2$ has an $tx_3$-hamiltonian path $P_1$ containing $tw\in E(G)$. If $G_0=\emptyset$, then we set $P=P_2\cup P_1\cup \{x_3x_1\}$. If $G_0=t'v_1$, then 
we set $P=P_2\cup P_1\cup \{x_3x_1, w'v_1, v_1t'\}-\{t'w'\}$. Otherwise $(G_0)^2$ has a hamiltonian cycle $C_0$ containing $t'w^*\in E(G)$ by 
Theorem \ref{fletheorem3} or Corollary \ref{flecor}(i), respectively, and we set $P=P_2\cup C_0\cup P_1\cup \{x_3x_1,w'w^*\}-\{t'w',t'w^*\}$. 
In all cases $P$ is an ${\cal F}_4$ $x_1x_2$-hamiltonian path in $G^2$ containing $x_3x_1, x_4w_4$ which are edges of $G$ as required. Note that if $t=x_4$, 
then $x_4w_4=tw$.

\vspace{2mm}
If $t'=x_4$, we proceed analogously as in the previous case with $G_1$ and $G_0$ switching roles.

\vspace{3mm} This completes the proof of Theorem \ref{dt}. \qed

\vspace{5mm}
\begin{center} 
 {\bf Acknowledgements} 
\end{center}

Research of the first author was supported by the FRGS Grant (FP036-2013B), the second author was supported by project P202/12/G061 of the Grant Agency of 
the Czech Republic, whereas research of the third author was supported in part by FWF-grant P27615-N25.

This publication was partly supported by the project LO1506 of the Czech Ministry of Education, Youth and Sports.

\end{document}